\documentclass[bj]{imsart}
 
\RequirePackage[OT1]{fontenc}
\RequirePackage{amsthm,amsmath}
\RequirePackage[colorlinks]{hyperref}
\usepackage{latexsym}
\usepackage{amsmath,amssymb,amsfonts}
 \usepackage[latin1]{inputenc}
\usepackage[T1]{fontenc}
\usepackage{amsmath,amssymb,amsfonts,color}
\def\vs#1{\vspace{#1 mm} }

\newtheorem{Theorem}{Theorem}[part]

\newtheorem{Proposition}{Proposition}[part]

\newtheorem{Lemma}{Lemma}[part]
\newtheorem{Corollary}{Corollary}[part]
\newtheorem{Remark}{Remark}[part]

\makeatletter \@addtoreset{equation}{section}

\@addtoreset{Definition}{section}

\@addtoreset{Theorem}{section}

\@addtoreset{Proposition}{section}

\@addtoreset{Property}{section}

\@addtoreset{Assumption}{section}

\@addtoreset{Corollary}{section}

\@addtoreset{Lemma}{section}

\@addtoreset{Remark}{section}

\@addtoreset{Example}{section}

 \def\vs#1{\vspace{#1mm}}

\def \E{\mathbb{E}}
\def \F{\mathbb{F}}

\def \M{\mathbb{M}}
\def \N{\mathbb{N}}
\def \R{\mathbb{R}}

\def\P{\mathbb{P}}
\def\Q{\mathbb{Q}}

\def\Xb{\bar X}

\def\Ac{{\cal A}}

\def\Cc{{\cal C}}
\def\Dc{{\cal D}}
\def\Ec{{\cal E}}
\def\Fc{{\cal F}}

\def\Hc{{\cal H}}

\def\Lc{{\cal L}}

\def\Nc{{\cal N}}
\def\Oc{{\cal O}}
\def\Pc{{\cal P}}

\def\Rc{{\cal R}}
\def\Sc{{\cal S}}

\def\Xc{{\cal X}}
\def\Yc{{\cal Y}}

\def\ep{\hbox{ }\hfill$\Box$}
\def\reff#1{{\rm(\ref{#1})}}
\def\be{\begin{eqnarray}}
\def\ee{\end{eqnarray}}
\def\beq{\begin{equation}}
\def\eeq{\end{equation}}
\def\b*{\begin{eqnarray*}}
\def\e*{\end{eqnarray*}}

\def\x{\times}

\def \sle{\;\le\;}
\def\={\;=\;}

\def\proof{{\noindent \bf Proof. }}

\def\And{\;\mbox{ and }\;}
\def\pourtout{\mbox{ for all } }
\def\Pas{\mathbb{P}-\mbox{a.s.}}
\def\no{\noindent}
\def\.{\;.}

\def\eps{\varepsilon}
\def\vp{\varphi}
\def\1{{\bf 1}}

\def\Esp#1{\mathbb{E}\left[#1\right]}

\def\Pro#1{\mathbb{P}\left[{#1}\right]}

\def\scap#1#2{\langle #1 , #2\rangle}

\def\Tr#1{\mbox{\rm Tr}\left[#1\right]}

\def\ti{{t_i}}
\def\tip{ {t_{i+1}} }

\def\dO{{\partial{\cal O}}}

\def\cD{{\mathcal{D}}}
\def\cK{{\mathcal{K}}}
\def\cC{{\mathcal{C}}}
\def\cO{{\mathcal{O}}}
\def\cA{{\mathcal{A}}}
\def\cN{{\mathcal{N}}}

\def\Xe{X^\varepsilon}
\def\Ye{Y^\varepsilon}

\def\Ze{Z^\varepsilon}

\def\te{\tau^\varepsilon}
\def\leftB{[\![}
\def\rightB{]\!]}

\def\eqref#1{\reff{#1}}

\usepackage{epsfig}
\usepackage{psfrag}

\def\cO{{\cal O}}

\def\Xb{\bar X}
\def\taub{{\bar \tau}}

\def\A#1{\textbf{(#1)}}
\def\Qro#1{\Q\left[#1\right]}
 \def\scap#1#2{ \langle #1, #2\rangle} 
\def\EspQ#1{\E^\Q\left[#1\right]}

\begin{document}
\begin{frontmatter}
\title{Strong Approximations of BSDEs in a domain}
\runtitle{Strong Approximations of BSDEs in a domain}

\begin{aug}
\author{\fnms{Bruno} \snm{Bouchard}\thanksref{a,e1} \ead[label=e1,mark]{bouchard@ceremade.dauphine.fr} \ead[label=u1,url]{www.ceremade.dauphine.fr/\~{}bouchard}}
\and 
\author{\fnms{Stéphane} \snm{Menozzi}\thanksref{b,e2}
\ead[label=e2,mark]{menozzi@math.jussieu.fr} 
\ead[label=u2,url]{www.proba.jussieu.fr/\~{}menozzi/}
}

\address[a]{CEREMADE,  Universit{\'e} Paris 9, place du Maréchal de Lattre de Tassigny, 75016 Paris France. \\ \printead{e1,u1}}
\address[b]{LPMA,  Universit{\'e}  Paris 7, 175 rue du Chevaleret, 75013 Paris France. \\
 \printead{e2,u2} }
 
\runauthor{B. Bouchard and S. Menozzi}

\affiliation{Some University and Another University}

\end{aug}

\begin{abstract}
We study the strong approximation of a Backward SDE with finite stopping time horizon, namely the first exit time of a forward SDE from a cylindrical domain. We use the Euler scheme approach of \cite{BoTo04,Zh04}. When the domain is piecewise smooth and under a non-characteristic boundary condition, we show that the associated strong error is at most of order $h^{\frac14-\eps}$ where $h$ denotes the time step and $\eps$ is any positive parameter. This rate corresponds to the strong exit time approximation.  It is improved to $h^{\frac12-\eps}$ when the exit time can be exactly simulated or  for a weaker form of the approximation error. Importantly, these results are obtained without uniform  ellipticity condition. 
\end{abstract}


\end{frontmatter}

\section{Introduction}

Let $T>0$ be a finite time horizon and $(\Omega,\Fc,\mathbb{P})$ be
a stochastic basis supporting a $d$-dimensional Brownian motion $W$.
We assume that the filtration $\F=(\Fc_t)_{t\le T}$ generated by $W$
satisfies the usual assumptions and that $\Fc_T = \Fc$.
 
Let $(X,Y,Z)$ be the solution of the decoupled Brownian Forward-Backward SDE 
	\be\label{eq def SDE intro}
	X_t&=& X_0 + \int_0^t b(X_s) ds + \int_0^t \sigma(X_s) dW_s 
	\\
	\label{eq def BSDE intro}
	Y_{t} &=& g(\tau,X_\tau) + \int_{t}^T \1_{s< \tau} f(X_s,Y_s,Z_s) ds -  \int_{t}^T Z_s dW_s\;\;\;,\;t\in [0,T]\;,
	\ee 
where $\tau$ is the first exit time of $(t,X_t)_{t\le T}$ from a cylindrical domain $D=[0,T)\x \Oc$ for some open piecewise smooth connected set $\Oc\subset \R^d$, and $b$, $\sigma$, $f$ and $g$ satisfy the usual Lipschitz continuity assumption. 

This kind of systems appears in many applications. In particular, it is well known that it is related to the solution of the semi-linear Cauchy Dirichlet problem 
	\be
	-\Lc u - f(\cdot,u,Du\sigma)= 0 \;\mbox{ on }\; D 
	&,& u= g \;\mbox{ on } \; \partial_p D\;,\label{eq edp intro}
	\ee
where $\Lc$ is the (parabolic) Dynkin operator associated to $X$,  i.e. for $\psi \in C^{1,2}$
	\b*
	\Lc \psi := \partial_t  \psi +	    \scap{b}{D \psi}  +  \frac12  \Tr{ a D^2\psi   } \;\;,\;\;a:=\sigma\sigma^*\;,
	\e*
and $\partial_p D:=\left([0,T)\x \partial \Oc\right) \cup \left(\{T\}\x \bar \Oc\right)$ is the parabolic boundary of $D$. More precisely, if the solution $u$ of \reff{eq edp intro} is smooth enough, then $Y=u(\cdot,X)$ and $Z=Du\sigma(\cdot,X)$. Thus, in the regular frame, solving \reff{eq def BSDE intro} is essentially equivalent to solving \reff{eq edp intro}. 

In this paper, we study an Euler scheme type approximation of \reff{eq def SDE intro}-\reff{eq def BSDE intro} similar to the one introduced in \cite{BoTo04,Zh04}, see also \cite{BoCh06,BoEl05,MaZh05}.  We first consider the Euler scheme approximation $\Xb$ of $X$ on some grid $\pi:=\{t_i=ih,\;i\le n\}$ with modulus $h:=T/n$, $n\in \N^*$. The exit time $\tau$ is approximated by the first discrete exit time   $\bar \tau$ of $(t_i,\Xb_\ti)_{t_i\in \pi}$  from $D$. Then, the backward Euler scheme of   $(Y,Z)$ is defined for $i=n-1,\ldots, 0$ as 
	\b*
	\bar Y_{\ti}
	 := 
	\Esp{\bar Y_\tip~|~\Fc_\ti}
	+
	\1_{\ti<\taub} \;h \;f(\bar X_{t_i},\bar Y_{t_i},\bar Z_{t_i})\label{eq def Euler BSDE intro}
	&,&
	\bar Z_{\ti} := h^{-1} \Esp{\bar Y_{\tip}\left(W_{\tip}-W_\ti\right)~|~\Fc_\ti}\;,
	\e*
with the terminal condition
	$
	\bar Y_T = g(\taub,\bar X_\taub)\;. 
	$
Here, $g$ is a suitable extension of the boundary condition on the whole space $[0,T]\x \R^d$. 

The main purpose of this paper is to provide bounds for the (square of the)  discrete time approximation error up to a stopping time $\theta \le T$ $\Pas$   defined as
	\be
	\label{eq def Err}
	{\rm Err}(h)_\theta^2
	:=
	\max_{i<n} \Esp{\sup_{t\in [\ti,\tip]} \1_{t\le \theta} |Y_t-\bar Y_\ti|^2}+ \Esp{\int_0^\theta \|Z_t-\bar Z_{\phi(t)}\|^2dt}\;,
	\ee
where $\phi(t):=\sup\{s\in\pi: s\le t\} $.

We are interested in  two important   cases:   $\theta=T$ and $\theta=\tau\wedge \taub$. 
The quantity 	${\rm Err}(h)_T$ coincides with the usual strong approximation error computed up to $T$.   The term ${\rm Err}(h)_{\tau\wedge \taub}$ should be more considered as a weak approximation error, since the length of the random time  interval $[0,\tau\wedge \taub]$ cannot be controlled sharply in pratice.  It essentially provides a bound for $Y_0-\bar Y_0$, or equivalently in terms of \reff{eq edp intro}, $u(0,X_0)-\bar Y_0$. Let us mention that a precise analysis of the weak error has been carried out by Gobet and Labart in \cite{GoLa07}  in the uniformly elliptic case with  $\Oc=\R^d$. 

\vs2

As in \cite{BoTo04}, \cite{MaZh02} and \cite{Zh04}, who also considered the limit case $\Oc=\R^d$ (i.e. $\tau=T$), 
the approximation  error can be naturally related to the error due to the approximation of $X$ by $\bar X_{\phi}$ and 
the regularity of the solution $(Y,Z)$ of \reff{eq def BSDE intro} through the quantities:
	\b*
 	\Rc(Y)_{\Sc^2}^{\pi}:= \max_{i<n}\Esp{\sup_{t\in [\ti,\tip]} |Y_{t}-Y_{\ti}|^2} &\And& \Rc(Z)_{\Hc^2}^{\pi}:= \Esp{\int_0^T \|Z_{t}-\hat Z_{\phi(t)}\|^2dt}
 	\e*
where 	
	\be\label{eq def hat Z}
	\hat Z_{\ti}:=h^{-1}\Esp{\int_\ti^\tip Z_s ds~|~\Fc_\ti} \;\mbox{ for }\; i<n\;.
	\ee 
In the case $f=0$, $Y$ is a martingale and $Y_\ti$ is the best $L^2$ approximation of $Y_t$ on the time interval $[\ti,\tip]$ by an $\Fc_\ti$-measurable random variable. In this case, Doob's inequalities imply that
$\Esp{\sup_{t\in [\ti,\tip]}|Y_t-\bar Y_\ti|^2}\ge \Esp{|Y_{\tip}-Y_{\ti}|^2}\ge c\;\Esp{\sup_{t\in [\ti,\tip]} |Y_{t}-Y_{\ti}|^2}$, for some universal constant $c>0$.  

Moreover, the definition \reff{eq def hat Z} implies that $\hat Z_\phi$ is the best approximation in $L^2([0,T]\x \Omega,dt\otimes d\P)$ of $Z$ by a process which is constant on each time interval $[\ti,\tip)$. Thus, $\Rc(Z)_{\Hc^2}^{\pi}$ $\le$ $\Esp{\int_0^T \|Z_t-\bar Z_{\phi(t)}\|^2dt}$. 

This justifies why $\Rc(Y)_{\Sc^2}^{\pi}$ and $\Rc(Z)_{\Hc^2}^{\pi}$ should play a crucial role in the convergence rate of ${\rm Err}(h)$ to $0$ as $h\to 0$.

Bounds for similar quantities   have previously been studied in \cite{BoTo04,Zh04} in the case $\Oc=\R^d$ and in \cite{BoCh06,MaZh05} in the case of reflected BSDEs. All these articles use a Malliavin calculus approach to derive a particular representation of $Z$. Due to the exit time, these techniques fail
in our setting. We propose a different approach that  relies on mixed analytic/probabilistic arguments.
Namely, we first adapt some barrier techniques from the PDE literature, see e.g. Chapter 14 in \cite{GiTr98} and Section \ref{sec gradient bound} below, to provide a bound for the modulus of continuity of $u$ on the boundary,   and  then some stochastic flows and martingale arguments to obtain 
an interior control on this modulus. Under the standing assumptions of Section \ref{sec notations and assumptions}, it allows to derive  that $\Rc(Y)_{\Sc^2}^{\pi}$ $+$ $\Rc(Z)_{\Hc^2}^{\pi}$ $=$ $O(h)$ and that $u$ is $1/2$-H${\rm \ddot{o}}$lder in time and Lipschitz continuous in space. 

\vs2 

To derive our final error bound on Err$(h)_\theta$, we additionally have to take into consideration the error 
coming from the approximation of $\tau$ by $\taub$. We show that $\Esp{|\tau-\taub|}=O(h^{\frac12-\eps})$ for all $\eps>0$. Combined with the previous controls on $\Rc(Y)^\pi_{\Sc^2}$ and $\Rc(Z)^\pi_{\Hc^2}$, this allows us to show that   ${\rm Err}(h)_{T}=O(h^{\frac14-\eps})$. Exploiting an additional control on a weaker form of error on $\tau-\taub$, we also derive that  ${\rm Err}(h)_{\tau\wedge \taub}=O(h^{\frac12-\eps})$. 
As a matter of facts, the global error is driven by the approximation error of the exit time which propagates backward thanks to the Lipschitz continuity of $u$.  

\vs2

Importantly, we do not assume specific non degeneracies of the diffusion coefficient but only a uniform non characteristic boundary condition and  uniform ellipticity close to the corners, recall that $\Oc$ is piecewise smooth. Using the transformation  proposed in \cite{Ko00}, these results could be extended to drivers with quadratic growth (for a bounded boundary condition $g$). Also, without major difficulties,
our results could be extended to time dependent domains and coefficients ($b$, $\sigma$ and $f$) under natural assumptions on the time regularity. We restrict here to the homogeneous cylindrical case for simplicity. 
\vs2

We note that the numerical implementation of the above scheme requires the approximation of the involved   conditional expectations. It can be performed by non-parametric regression techniques, see  e.g. \cite{GoLeWa05} and \cite{LoSc01}, or a quantization approach, see e.g.  \cite{BaPa02}   and \cite{DeMe06,DeMe07}.  In both cases, the additional error is analyzed  in the above papers and can be extended to our framework. We note that the Malliavin approach of \cite{BoTo04} cannot be directly applied here due to the presence of the exit time. Concerning a direct computable algorithm, we mention the work of Milstein and Tretyakov \cite{MiTr01} who use  a simple random walk approximation of the  Brownian motion.  However, their results require  strong smoothness assumptions on the solution of \reff{eq edp intro} as well as a uniform ellipticity condition. 

\vs2

The rest of the paper is organized as follows. We start with some notations and assumptions in Section \ref{sec notations and assumptions}. Our main results are presented in Section \ref{sec main results}.
In Section \ref{sec prop majo erreur en terme de regu}, we provide a first bound on the error: it involves the error due to the discrete time approximation of $\tau$ by $\taub$ and the regularity of the solution $(Y,Z)$ of \reff{eq def BSDE intro}. The discrete approximation of $\tau$ is specifically studied in  Section \ref{sec thm main erreur ta}.  Eventually,  Section \ref{sec regu} is devoted to the analysis of the regularity of \reff{eq edp intro} and \reff{eq def BSDE intro} under our current assumptions.

\section{Notations and assumptions}\label{sec notations and assumptions}

Any element $x \in \R^d$, $d\ge 1$,  will be
identified to a line vector with $i$-th component $x^i$ and
Euclidean norm $\|x\|$.   The
scalar product on $\R^d$ is denoted by $\scap{x}{y}$. The open ball of center $x$ and radius $r$ is denoted by $B(x,r)$, $\bar B(x,r)$ is its closure. 
Given a non-empty set $A \subset \R^d$, we similarly denote by  $B(A,r)$ and $\bar B(A,r)$ the sets $\{x\in \R^d~:~d(x,A)< r\}$ and $\{x\in \R^d~:~d(x,A)\le  r\}$ where $d(x,A)$ stands for the Euclidean distance of $x$ to $A$. 
 For a $(m\x d)$-dimensional matrix $M$, we denote $M^*$ its transpose and we write $M\in \M^d$ if
$m=d$.   For a smooth function $f(t,x)$,   $Df$ and $D^2 f$ stand for its gradient (as a line vector) and Hessian matrix with 
respect to its second component. If it depends on some extra components, we denote by  $\partial_t f(t,x,y,z)$, $\partial_x f(t,x,y,z)$, etc... its partial gradients.

\subsection{Euler scheme approximation of BSDEs}\label{sec Euler scheme approximation of BSDEs}
 
From now on, we assume that the coefficients of  \reff{eq def SDE intro}-\reff{eq def BSDE intro} satisfy: 

\vs2
 
\no \A{HL}: There is a constant $L>0$ such that for all $(t,x,y,z,t',x',y',z') \in ([0,T]\x\R^d\x \R\x \R^d)^2$:
	\b*
	\left\|(b,\sigma,g,f)(t,x,y,z)-(b,\sigma,g,f)(t',x',y',z')\right\| &\le& L\; \left\|(t,x,y,z)-(t',x',y',z')\right\|\;,
	\\
	\left\|(b,\sigma,g,f)(t,x,y,z)\right\| &\le&  L\; \left(1+\left\|(x,y,z)\right\|\right)\;.
	\e*

Under this assumption, it is well known, see e.g. \cite{PaPe92,Pe91}, that we have existence and uniqueness of a solution $(X,Y,Z)$ in $\Sc^2\x\Sc^2\x \Hc^2$, where   we  denote by $\Sc^2$  the set of real valued adapted continuous processes $\xi$ satisfying
    $
    \|\xi\|_{\Sc^2} :=  \Esp{\sup_{ t\le T } |\xi_t|^2}^{\frac12} \;<\;\infty\;
    $, 
and by $\Hc^2$ the set of progressively measurable $\R^d$-valued processes $\zeta$ for which
    $
    \|\zeta\|_{\Hc^2} :=  \E[  \int_0^T |\zeta_t|^2 dt ]^{\frac12} \;<\;\infty\;.
    $

\vs2

As usual, we shall approximate the solution of \reff{eq def SDE intro} by its Euler scheme $\bar X$ associated to a grid 
	$$
	\pi:=\{t_i=ih\;,\;i\le n\}\;,\; h:=T/n\;,\;n\in \N^*\;,
	$$ 
defined by 
	\be\label{eq def Euler SDE}
	\bar X_t&=& X_0 + \int_0^t b(\bar X_{\phi(s)}) ds + \int_0^t \sigma(\bar X_{\phi(s)}) dW_s \;\;,\;t\ge 0\;,
	\ee
where we recall that 
	$ 
	\phi(s):={\rm arg}\max\{t_i,\;i\le n~:~t_i\le s\} \;\;\mbox{ for } s\ge 0\;. 
	$

Regarding the approximation of \reff{eq def BSDE intro}, we adapt the approach of \cite{Zh04} and \cite{BoTo04}. First, we approximate the exit time $\tau$ by the first exit time of the Euler Scheme $(t,\bar X_{t})_{t\in \pi}$ from $D$ on the grid $\pi$:
	\b*
	\taub:=\inf\{t\in \pi~:~\bar X_{t}\notin \Oc\}\wedge T\;. 
	\e* 

\begin{Remark}{\rm
Note that one could also approximate $\tau$ by  $\tilde \tau:=\inf\{t\in [0,T]~:~\bar X_{t}\notin \Oc\}\wedge T$, the first exit time of the ``continuous version'' of the Euler scheme  $(t,\bar X_{t})_{t\in [0,T]}$, as it is done for linear problems, i.e. $f$ is independent of $(Y,Z)$, see e.g. \cite{Go00}. However, in the case where $\Oc$ is not a half-space, this requires  additional local approximations of the boundary by tangent hyperplanes and will not allow to improve our strong approximation error, compare Corollaire 2.3.2. in \cite{Go98} with  Theorem \ref{thm main erreur ta} below. 
}
\end{Remark}

Then, we define the discrete time process $(\bar Y,\bar Z)$ on $\pi$  by 
	\be
	\bar Y_{\ti}
	&:=&
	\Esp{\bar Y_\tip~|~\Fc_\ti}
	+
	\1_{\ti<\taub} \;h \;f(\bar X_{t_i},\bar Y_{t_i},\bar Z_{t_i})\;,\label{eq def Euler BSDE}
	\\
	\bar Z_{\ti}&:=&h^{-1} \Esp{\bar Y_{\tip}\left(W_{\tip}-W_\ti\right)~|~\Fc_\ti} \;\;,\;i<n\;,
	\label{eq def Euler BSDE bar Z}
	\ee
with the terminal condition
	\be
	\bar Y_T&=&g(\taub,\bar X_\taub)\;. 
	\ee
Observe that $\bar Y_\ti\1_{\ti\ge \taub}=g(\taub,\bar X_\taub)\1_{\ti\ge \taub}$ and that $\bar Z_{\ti}\1_{\ti\ge \taub}=0$.

\vs2

One easily checks that  $(\bar Y_\ti,\bar Z_\ti)\in L^2$ for all $i\le n$ under \A{HL}. It then follows from the martingale representation theorem that we can find $\tilde Z \in \Hc^2$ such that 
	\be\label{eq def tilde Z}
	\bar Y_{\tip}-\Esp{\bar Y_\tip~|~\Fc_\ti}&=& \int_\ti^\tip \tilde Z_s dW_s \;\;\pourtout\;i<n\;. 
	\ee
This allows us to consider a continuous time extension of $\bar Y$ in $\Sc^2$ defined on $[0,T]$ by
	\be\label{eq def Euler BSDE avec tilde Z}
	\bar Y_t&=& g(\taub,\bar X_\taub) +\int_t^T \1_{s<\taub}\; f(\bar X_{\phi(s)},\bar Y_{\phi(s)},\bar Z_{\phi(s)}) ds - \int_t^T \tilde Z_s dW_s\;. 
	\ee

\begin{Remark}\label{rem Z=0}{\rm Observe that $Z=0$ on $]\tau,T]$ and $\tilde Z=0$ on $]\taub,T]$.  For later use, also notice that the Itô isometry and \reff{eq def tilde Z} imply
	\be\label{eq lien bar Z tilde Z}
	\bar Z_\ti&=& h^{-1}\; \Esp{\int_\ti^\tip \tilde Z_s ds ~|~\Fc_\ti} \;\;,\;i<n\;. 
	\ee
}
\end{Remark}

\subsection{Assumptions on $\Oc$, $\sigma$ and $g$} \label{sec Assumptions}

Our main result holds under some additional assumptions on $\Oc$, $\sigma$ and $g$. Without loss of generality, we can specify them in terms of the constant $L$ which appears in \A{HL}. 

\vs2

We first assume that the domain $\Oc$ is a finite intersection  of smooth domains with compact boundaries:

\vs2

\no \A{D1}: We have $\Oc:=\bigcap_{\ell=1}^m  \Oc^\ell$ where $m\in \N^*$ and   $\Oc^\ell$ is a $C^2$
 domain of $\R^d$ for each $1\le \ell\le m$. Moreover, $\Oc^\ell$ has a compact boundary, $\sup\{\|x\|~:~x\in \partial \Oc^\ell\}\le L$, for each $1\le \ell \le m$.

\vs2

It follows from  Appendix 14.6 in \cite{GiTr98} that there is a function $d$ which coincides with the algebraic distance to $\partial \Oc$, in particular
	$\Oc:=\{x\in \R^d~:~d(x)>0\}\;$, 
and is $C^2$ outside of a neighborhood  $B(\Cc,L^{-1})$ of the set of corners
	\b*
	\Cc:=\bigcap_{\ell\ne k=1}^m  \partial \Oc^\ell\cap \partial \Oc^k\;.
	\e*  
 
\vs2

We also assume that the domain satisfies a uniform exterior sphere condition as well as a uniform truncated interior cone condition:

\vs2

\no {\A{D2}:} For all $x\in \partial \Oc$, there is $y(x)\in  \Oc^c$, $r(x)\in [L^{-1},L]$ and $\delta(x)\in B(0,1)$  such that 
	\b*
 	\bar B(y(x),r(x))\cap \bar \Oc=\{x\} \;
	\mbox{ and } \{x'\in    B(x,L^{-1})~:~\scap{x'-x}{\delta(x) }\ge (1- L^{-1})\|x'-x\|\}\subset \bar \Oc\;.
	\e*

In view of \A{D1}, these last assumptions are actually automatically satisfied outside a neighborhood of the set of corners, see  e.g. Appendix 14.6 in \cite{GiTr98}.

\vs2

In order to ensure that the associated first boundary value problem is well posed in the (unconstrained) viscosity sense, we shall also assume that 
	\b*
	a:=\sigma\sigma^*
	\e*
satisfies a non-characteristic boundary condition outside the set of corners $\Cc$ and a uniform ellipticity condition on a neighborhood of $\Cc$:

\vs2

\no {\A{C}:} We have 
	\b*
	\inf\{n(x)a(x) n(x)^*~:~x \in \partial \Oc\setminus  B(\Cc,L^{-1})  \} \ge L^{-1}\; \mbox{ where } n(x):=Dd(x)\;,
	\e*
and	
	\b*
	\inf\{\xi a(x)\xi^*~:~\xi \in \partial B(0,1)\;,\;x \in \bar \Oc\cap B(\Cc,L^{-1}) \} \ge L^{-1}\;.
	\e*

In particular, it guarantees that the process $X$ is non-adherent to the boundary.

Observe that $n$ coincides with the inner normal unit on $\partial \Oc$ outside the set of corners. By abuse of notations, we write  $n(x)$ for $Dd(x)$, whenever this quantity  is well defined,  even if $x$ $\notin \partial \Oc$. 

\vs2

Importantly, we do not assume that $\sigma$ is non degenerate in the whole domain. 

\vs2

We finally assume that $g$ is smooth enough:

\vs2

\no \A{Hg}: $g \in C^{1,2}([0,T]\x\R^d)$ and 
	$
	\|\partial_t g\|+\|Dg\|+\|D^2g\| \sle L \;\;\;\mbox{ on }  [0,T]\times \R^d\;. 
	$
\\

Clearly, this smoothness assumption could be imposed only on a neighborhood of  $\partial \Oc$. Since it is compact and $Y$ depends  on $g$ only on $\partial \Oc$, we can always construct a suitable extension of $g$ on $\R^d$ which satisfies the above condition. Actually, one could only assume that $g$ is Lipschitz in $(t,x)$ and has a Lipschitz continuous derivative in $x$. With this slightly weaker condition, all our arguments would go through after possibly replacing $g$ by a sequence of regularized versions and then passing to the limit, see  Section  \ref{subsec regularization procedure} for similar kind of arguments. 

\section{Main results}\label{sec main results}
We first provide a general   control on the quantities in \reff{eq def Err} in terms of $\Rc(Y)_{\Sc^2}^\pi$, $\Rc(Z)_{\Hc^2}^\pi$ and  $|\tau-\taub|$. Let us mention that this type of result is now rather standard when $\Oc=\R^d$, see e.g. \cite{BoTo04}, and requires only the Lipschitz continuity assumptions of \A{HL}. 
 
\begin{Proposition}\label{prop majo erreur en terme de regu} Assume that \A{HL} and \A{Hg}   hold. Then, there exist $C_L>0$ and a positive random variable  $\xi_L$ satisfying $\Esp{(\xi_L)^p}\le C_L^p$ for all $p\ge 2$ such that 
	\be
	{\rm Err}(h)_T^2
	 \sle  
	C_L\left(h + \Rc(Y)_{\Sc^2}^\pi +\Rc(Z)_{\Hc^2}^\pi +\Esp{\xi_L|\tau-\taub|+\1_{\taub<\tau} \int_{\taub }^{\tau}   \|Z_s\|^2  ds} \right)
	\; \label{eq majo erreur en terme de regu T}
	\ee
and 
	\be
	{\rm Err}(h)_{\tau\wedge \taub}^2 \le {\rm Err}(h)_{\tau_+\wedge \taub}^2 
	&\le& 
	C_L\left(h + \Rc(Y)_{\Sc^2}^\pi +\Rc(Z)_{\Hc^2}^\pi \right)+ \Esp{ \Esp{\xi_L|\tau-\taub|~|~\Fc_{\tau_+\wedge \taub}}^2}
	\nonumber\\
	&+&
	C_L\;
	\Esp{\1_{\taub<\tau}\Esp{ \int_\taub^\tau   \|Z_s\|   ds~|~\Fc_{\taub}}^2}
	\;. \label{eq majo erreur en terme de regu tau+ taub}
	\ee
where $\tau_+$ is the next time after $\tau$ in the grid $\pi$:
	$
	\tau_+:=\inf\{t\in \pi~:~\tau\le t\}\;. 
	$	
\end{Proposition} 
The proof will be provided in Section \ref{sec prop majo erreur en terme de regu} below.  Note that we shall control ${\rm Err}(h)_{\tau\wedge \taub}^2$ through the slightly stronger term ${\rm Err}(h)_{\tau_+\wedge \taub}^2$, see \reff{eq majo erreur en terme de regu tau+ taub}. This will allow us to work with stopping times with values in the grid $\pi$ which will be technically easier, see Remark \ref{rem dble int} below.

\vs2
 
In order to provide a convergence rate for ${\rm Err}(h)_T^2$ and ${\rm Err}(h)_{\tau_+\wedge \taub}^2$, it remains to control the quantities  $\Rc(Y)_{\Sc^2}^\pi$, $\Rc(Z)_{\Hc^2}^\pi$ and the terms involving the difference between $\tau$ and $\taub$. 

\vs2
The error due to the approximation of $\tau$ by $\taub$ is controlled by the following estimate that extends to the non uniformly elliptic case previous results obtained in \cite{Go98}, see its Corollaire 2.3.2. The proof of this Theorem is provided in Section \ref{sec thm main erreur ta} below. 

\begin{Theorem}\label{thm main erreur ta} Assume that $b$ and $\sigma$ satisfy \A{HL} and that \A{D1} and \A{C}  hold. Then, for $\eps\in (0,1)$  and each positive random variable $\xi$ satisfying $\Esp{(\xi)^p}\le C_L^p$ for all $p\ge 1$, there is $C_L^\eps>0$ such that
	\b*
	\Esp{\Esp{\xi\; |\tau-\bar \tau|~|~\Fc_{\tau_+\wedge \taub} }^2} &\le& C_L^\eps h^{1-\eps}\;. 
	\e*
In particular, for each $\eps\in (0,1/2)$, there is $C_L^\eps>0$ such that 
	\b*
	\Esp{|\tau-\bar \tau|} &\le& C_L^\eps h^{1/2-\eps}\;. 
	\e*
\end{Theorem}

In  \cite{Go98}, the last bound is derived under   a uniform ellipticity condition on $\sigma$ and cannot be exploited in our setting, recall that we only assume \A{C}. Up to the $\eps$ term, it can not be improved. Indeed, in the special case of a uniformly elliptic diffusion in a smooth bounded domain, it has been shown in \cite{GoMe07} that $\Esp{\tau -\bar \tau}=C h^\frac12+o(h^\frac12)$ for some $C>0$, see Theorem 2.3 of this reference. 
\\

Our next result concerns the regularity of $(Y,Z)$ and is an extension to our framework of similar results obtained in \cite{MaZh02}, \cite{BoTo04}, \cite{BoEl05} and \cite{BoCh06} in different contexts.

\begin{Theorem}\label{thm main regu} Let the conditions \A{HL}, \A{D1}, \A{D2}, \A{C}  and \A{Hg} hold.  
Then, 
	\be\label{eq main thm regu}
  \Rc(Y)_{\Sc^2}^\pi+ \Rc(Z)_{\Hc^2}^\pi
 	&\le&
	C_L\;h\;.
	\ee
Moreover, for all  stopping times $\theta,\vartheta$ satisfying $\theta\le \vartheta\le T$ $\Pas$, one has 
	\be\label{eq main thm regu 3}
	\Esp{\sup_{\theta\le s \le \vartheta} |Y_s-Y_\theta|^{2p} } &\le& \Esp{\xi_L^p\;|\vartheta-\theta|^p}\;,\;\;p\ge 1\;,
	\ee
and  
	\be\label{eq main thm regu 2}
	\Esp{\int_\theta^\vartheta \|Z_s\|^p ds~|~\Fc_\theta} &\le& \Esp{\xi^p_L|\vartheta-\theta|~|~\Fc_\theta}\;\;,\;p=1,2\;,
	\ee
where  $\xi_L^p$ is a positive random variable which satisfies
	$
	\Esp{|\xi_L^p|^q} <\infty\;,
	$
for all $q\ge 1$.  

In addition, the unique continuous viscosity solution $u$ of 
\reff{eq edp intro},  in the class of continuous solutions with polynomial growth, is uniformly $1/2 $-H\"older continuous in time and Lipschitz continuous in space, i.e.
	\be\label{eq regu u}
	|u(t,x)-u(t',x')|\le C_L\left(|t-t'|^\frac12 + \|x-x'\|\right)\;\;\;\pourtout\; (t,x)\And(t',x')\in \bar D\;. 
	\ee 
\end{Theorem}
The proof is provided in Section \ref{sec regu} below.  The bound \reff{eq main thm regu 2} can be interpreted as a weak bound on the gradient, whenever it is well defined, of the viscosity solution of \reff{eq edp intro}. It implies that $Y$ is $1/2$-Hölder continuous in $L^2$ norm. This result is rather standard under our Lipschitz continuity assumption in the case where $\Oc=\R^d$, i.e. $\tau=T$, but seems to be new in our context and under our assumptions.  The bound $\Rc(Z)_{\Hc^2}^\pi \le C_L\;h$ can be seen as a weak regularity result on this gradient. It would be straightforward if one could show  that $Du\sigma$ is uniformly  $1/2$-Hölder in time and Lipschitz in space, which is not true in general.

\vs2 

Combining the above estimates, we finally obtain our main result which provides an upper bound for the  convergence rate of ${\rm Err}(h)_{\tau_+\wedge \taub}^2$ (and thus for ${\rm Err}(h)_{\tau\wedge \taub}^2$) and ${\rm Err}(h)_{T}^2$. 

\begin{Theorem} Let the conditions \A{HL}, \A{D1}, \A{D2}, \A{C}  and \A{Hg} hold.  
	Then, for each $\eps \in (0,\frac12)$, there is $C^\eps_L>0$ such that 
	\b* 
	{\rm Err}(h)_{\tau_+\wedge \taub}^2 \sle C^\eps_L\;h^{1-\eps}
	&\And&
	{\rm Err}(h)_{T}^2 \sle C^\eps_L\;h^{\frac12-\eps}\;. 
	\e*
\end{Theorem}

This extends   the results of \cite{BoCh06,BoEl05,Zh04} who obtained similar bounds in different contexts.

\begin{Remark}{\rm When $\tau$ can be exactly simulated, we can replace $\bar \tau$ by $\tau$ in the scheme \reff{eq def Euler BSDE}-\reff{eq def Euler BSDE bar Z}. 
In this case, the two last terms in the right hand-sides of \reff{eq majo erreur en terme de regu T} and \reff{eq majo erreur en terme de regu tau+ taub} cancel and we retrieve the convergence rate of the case $\Oc=\R^d$, see e.g. \cite{BoTo04}. }
\end{Remark} 

\begin{Remark}{\rm
 Note that the Lipschitz continuity assumption with respect to the $x$ variable on $g$ and $f$ is  only used to control at the right order the error term coming from the approximation of $X$ by $\bar X$ in $g$ and $f$.  
 If one is only interested in the convergence of ${\rm {E}rr}(h)_{T}$ this assumption can be weakened. Indeed, if we only assume that 

\vs2
\A{HL'$_1$}:  $b$, $\sigma$ satisfy \A{HL},  $\sup\{|f(\cdot,y,z)|,$ $(y,z)\in \R\x \R^d\}$ and $g$ have polynomial growth, and   $f(x,\cdot)$ is uniformly Lipschitz continuous, uniformly in $x\in \R^d$,  
\vs2

  a weak version of \reff{eq majo erreur en terme de regu T} can still be established up to an obvious modification of the proof of Proposition \ref{prop erreur BSDE theta} below. Namely,   there exists $C>0$ and a positive random variable  $\xi$ satisfying $\Esp{(\xi)^p}\le C_L^p$ for all $p\ge 2$ for which
\be
	{\rm Err}(h)_T^2
	 &\sle&  
	C\left(h +   \Esp{\int_0^T |Y_s-Y_{\phi(s)}|^2 ds} +\Rc(Z)_{\Hc^2}^\pi+ \Esp{ \xi |\tau-\taub| +\int_0^T\1_{\taub<\tau} \int_{\taub }^{\tau}   \|Z_s\|^2  ds} \right)~~~~~\nonumber	\\
	&+&C \Esp{|g(\tau,X_\tau)-g(\taub,\bar X_\taub)|^2 +\int_0^T |f(X_s,Y_s,Z_s)-f(\bar X_{\phi(s)},Y_s,Z_s)|^2  ds} \;.\label{eq majo erreur en terme de regu T bis}
	\ee
The terms  $\Esp{\int_0^T |Y_s-Y_{\phi(s)}|^2 ds}$ and $\Rc(Z)_{\Hc^2}^\pi$ are   easily seen to go $0$ with $h$, see e.g. the proof of Proposition 2.1 in \cite{BoEl05} for details. As for the other terms in the first line, it suffices to appeal to  Theorem \ref{thm main erreur ta} which implies that $ \Esp{\xi |\tau-\taub|}\to 0$ and that    $\bar \tau\to \tau$ in probability under \A{D1} and \A{C}. Note that the last assertion implies that  $\Esp{\int_0^T\1_{\taub<\tau} \int_{\taub }^{\tau}   \|Z_s\|^2  ds}\to 0$ and    $X_\tau-\bar X_\taub\to 0$  in probability. Hence, under the   additional continuity assumption

 \vs2
\A{HL'$_2$}: $g$ and $f(\cdot,y,z)$ are continuous, uniformly in $(y,z)\in  \R\x \R^d$,
\vs2

 we deduce that the two last terms in the second line go to $0$ as well.  
}
\end{Remark}


\section{Euler scheme approximation error: Proof of Proposition \ref{prop majo erreur en terme de regu}}\label{sec prop majo erreur en terme de regu}

In this section, we provide the proof of Proposition \ref{prop majo erreur en terme de regu}. 
We first recall  some standard controls on $X$, $(Y,Z)$ and $\bar X$ which holds under \A{HL}. 

From now on, $C^\eta_L$ denotes a generic constant whose value may change from line to line but which depends only on $X_0$, $L$ and some extra parameter $\eta$ (we simply write $C_L$ if it depends only on $X_0$ and $L$). Similarly, $\xi^\eta_L$ denotes a generic non-negative random variable such that $\Esp{|\xi^\eta_L|^p}\le C_L^{\eta,p}$ for all $p\ge 1$ (we simply write $\xi_L$ if it does not depend on the extra parameter $\eta$).

\begin{Proposition}\label{prop X bar X} Let \A{HL} hold. Fix $p\ge 2$.  Let $\vartheta$  be a stopping time  with values in $[0,T]$. Then 
	\b*
	\Esp{\sup_{t\in [\vartheta,T]}  \|Y_t\|^p +\left(\int_\vartheta^T \|Z_t\|^2 dt\right)^\frac{p}{2}~|~\Fc_\vartheta}&\le& C_L^p(1+\|X_\vartheta\|^p)\;
	\e*
and 
	\b*
	\Esp{\sup_{t\in [\vartheta,T]} \left(\|X_t\|^p+\|\Xb_t\|^p\right)~|~\Fc_\vartheta}&\le& \xi_L^p\;.
	\e*	
Moreover,
	\b*
	&\max\limits_{i<n}\Esp{\sup\limits_{t\in [\ti,\tip]} \left(\|X_t-X_\ti\|^p+\|\Xb_t-\Xb_\ti\|^p\right)}
	+ 
	\Esp{\sup\limits_{t\in [0,T]} \|X_t-\Xb_t\|^p}
	\le C_L^p h^\frac{p}{2}\;,&
	\\
	&\Pro{\sup\limits_{t\le T} \|\bar X_t-\bar X_{\phi(t)}\|> r} \le C_L\;r^{-4}\; h \;\;,\;r>0\;,&
	\e*
and, if  $\theta$ is a stopping time with values in $[0,T]$ such that $\vartheta\le \theta\le \vartheta+h$ $\Pas$, then
	\b*
	\Esp{ \|\Xb_\theta-\Xb_\vartheta\|^p+\|X_\theta-X_\vartheta\|^p~|~\Fc_\vartheta} 	&\le& \xi_L^p h^\frac{p}{2}\;. 
	\e*
\end{Proposition}

\begin{Remark}\label{rem hp sp bsde discret}{\rm For later use, observe that the Lipschitz continuity assumptions   \A{HL}   ensure  that 
	\b*
	\Esp{\sup_{t\in [\vartheta,T]}  \|\bar Y_t\|^p +\left(\int_\vartheta^T \|\tilde Z_t\|^2 dt\right)^\frac{p}{2}~|~\Fc_\vartheta}&<& \infty\;
	\;\;\mbox{ for all $p\ge 2$ .}
	\e* 
}
\end{Remark}

In order to avoid the repetition of similar arguments depending whether we consider ${\rm {E}rr}(h)^2_{\theta}$ with $\theta=T$ or $\theta=\tau_+\wedge \taub$, we first state an abstract version of  Proposition \ref{prop majo erreur en terme de regu} for some stopping time $\theta$ with values in $\pi$. 

\begin{Proposition}\label{prop erreur BSDE theta} Assume that $b,\sigma$ and $f$ satisfy \A{HL}. Then, for all stopping time $\theta$ with values in $\pi$, we have 
	\b*
		{\rm {E}rr}(h)^2_{\theta}
		\le
		C_L\left(
		h+ \Esp{ |Y_\theta-\bar Y_\theta|^2}
		+
		\Rc(Y)_{\Sc^2}^\pi+ \Rc(Z)_{\Hc^2}^\pi
		+ 
		\Esp{  \int_{\taub\wedge\tau \wedge \theta}^{(\taub\vee \tau) \wedge \theta}  \left( \xi_L  + \1_{\taub<\tau} \|Z_s\|^2 \right)ds}
		\right)
		\;.
	\e*
\end{Proposition} 

Let us first make the following Remark which will be of important use below. 

\begin{Remark}\label{rem dble int}{\rm Let $\vartheta\le \theta$ $\Pas$ be two stopping times with values in $\pi$ and $H$ be some adapted process in $\Sc^2$. Then, recalling that $\tip-\ti=h$, it follows from \reff{eq lien bar Z tilde Z} and Jensen's inequality  that 
	\b*
	\Esp{\int_\vartheta^\theta H_{\phi(s)} \|\bar Z_{\phi(s)}\|^2  ds }
	&=&
	\sum_{i<n} 
	\Esp{\int_\ti^\tip \1_{\vartheta\le \ti<\theta} \;  H_{\ti} \left\|\Esp{h^{-1} \;\int_\ti^\tip \tilde Z_{u}du~|~\Fc_\ti}\right\|^2 ds }
	\\
	&\le&
	\sum_{i<n} 
	\Esp{\int_\ti^\tip \1_{\vartheta\le \ti<\theta}  \; H_{\ti} h^{-1} \int_\ti^\tip \|\tilde Z_{u}\|^2 du ds }
	\\
	&\le&
	\Esp{\int_\vartheta^\theta H_{\phi(s)} \|\tilde Z_s \|^2  ds }\;. 
	\e*
By definition of $\hat Z$, see \reff{eq def hat Z}, the same inequality holds with $(\hat Z,Z)$ or $(\hat Z-\bar Z,Z-\tilde Z)$ in place of $(\bar Z,\tilde Z)$. This remark will allow us to control $\|Z-\bar Z_{\phi}\|$ through $\|Z-\tilde Z\|$ and  $\|Z-\hat Z_{\phi}\|$, see \reff{eq cont Z bar Zphi par Z hat Z et Z tilde Z} below, which is a key argument in the proof of Proposition \ref{prop erreur BSDE theta}. 
Observe that the above inequality does not apply if $\vartheta$ and $\theta$ do not take   values in $\pi$. This explains why it is easier to work with $\tau_+$ instead of $\tau$, i.e. work on ${\rm {E}rr}(h)^2_{\tau_+\wedge \taub}$ instead of ${\rm {E}rr}(h)^2_{\tau\wedge \taub}$. 
}
\end{Remark}

\no {\bf Proof of Proposition \ref{prop erreur BSDE theta}.} We adapt the arguments used in the proof of  Theorem 3.1 in \cite{BoTo04} to our setting. 
By applying Itô's Lemma to $(Y-\bar Y)^2$ on $[t\wedge \theta,\tip\wedge \theta]$ for $t\in [\ti,\tip]$ and $i<n$, we first deduce from \reff{eq def BSDE intro} and \reff{eq def Euler BSDE avec tilde Z} that 
	\b*
	\Delta^\theta_{t,\tip}&:=&
	\Esp{ |Y_{t\wedge \theta}-\bar Y_{t\wedge \theta}|^2 + \int_{t\wedge \theta}^{\tip\wedge \theta} \|Z_s-\tilde Z_s\|^2 ds }
	\\
	&=& 
	\Esp{ |Y_{\tip\wedge\theta}-\bar Y_{\tip\wedge\theta}|^2 }
	+
	\Esp{ 2 \int_{t\wedge \theta}^{\tip\wedge \theta} (Y_s-\bar Y_s) \left(\1_{s<\tau} f(\Theta_s)-\1_{s<\taub} f(\bar \Theta_{\phi(s)})\right)ds } 
	 \;, 
	\e*
where the martingale terms cancel thanks to Proposition \ref{prop X bar X} and Remark \ref{rem hp sp bsde discret}, and where $\Theta:=(X,Y,Z)$ and $\bar \Theta:=(\bar X,\bar Y,\bar Z)$. Using the inequality $2ab\le a^2+b^2$, we then deduce that, for $\alpha>0$  to be chosen later on,
	\b*
	\Delta^\theta_{t,\tip}
	&\le& 
	\Esp{ |Y_{\tip\wedge\theta}-\bar Y_{\tip\wedge\theta}|^2 }
	+
	 \alpha\; \Esp{  \int_{t\wedge \theta}^{\tip\wedge \theta}  |Y_s-\bar Y_s|^2ds }
	\nonumber\\
	&+& 2 \alpha^{-1}\Esp{ \int_{t\wedge \theta}^{\tip\wedge \theta} \1_{s<\taub} \left( f(\Theta_s)- f(\bar \Theta_{\phi(s)})\right)^2 ds 
	+
	\int_{t\wedge \theta}^{\tip\wedge \theta}\1_{\taub \le s<\tau}  \left( f(\Theta_s)\right)^2 ds } 
	\nonumber\\
	&+&
	2\alpha^{-1}\Esp{\int_{t\wedge \theta}^{\tip\wedge \theta}\1_{\tau \le s<\taub}  \left( f(\Theta_s)\right)^2 ds } \;. 
	\nonumber
	\e*
Recall from Remark \ref{rem Z=0} that $Z=0$ on $]\tau,T]$. Since $Y_t=g(\tau,X_\tau)$ on $\{t\ge \tau\}$, we then deduce from  
\A{HL} and  Proposition \ref{prop X bar X} that
	\be
	\Delta^\theta_{t,\tip}
	&\le&
	\Esp{ |Y_{\tip\wedge\theta}-\bar Y_{\tip\wedge\theta}|^2 }
	+
	 \alpha\;\Esp{  \int_{t\wedge \theta}^{\tip\wedge \theta}     |Y_s-\bar Y_s|^2ds } 
	\nonumber\\
	&+& 
	C_L\;\alpha^{-1}\Esp{  h\;|Y_{\ti\wedge \theta}-\bar Y_{\ti\wedge \theta}|^2  + \int_{t\wedge \theta}^{\tip\wedge \theta}     |Y_s- Y_{\phi(s)}|^2ds } 
	\nonumber\\
	&+& 
	C_L \;\alpha^{-1}
	\Esp{
	\int_{t\wedge \theta}^{\tip\wedge \theta}\left(h+ \|Z_s-\hat Z_{\phi(s)}\|^2+\|\hat Z_{\phi(s)}-\bar Z_{\phi(s)}\|^2 \right)  ds 
	} 
	\nonumber\\
	&+&  C_L\; \alpha^{-1} \Esp{ \int_{t\wedge \theta}^{\tip\wedge \theta} (\xi_L\1_{\tau\wedge \taub\le s\le \tau\vee\taub}  + \1_{\taub \le s<\tau}   \|Z_s\|^2 ) ds } \;. 
	\label{eq borne Delta 1}
	\ee
It then follows from Gronwall's Lemma that   
	\be
	\Esp{ |Y_{t\wedge \theta}-\bar Y_{t\wedge \theta}|^2 } 
	&\le& 
	(1+C^\alpha_L\;h) \Esp{ |Y_{\tip\wedge\theta}-\bar Y_{\tip\wedge\theta}|^2 } \nonumber\\
	&+& 
	(C_L\;\alpha^{-1}+C^\alpha_L\;h)
	\Esp{  h\;|Y_{\ti\wedge \theta}-\bar Y_{\ti\wedge \theta}|^2  + \int_{t\wedge \theta}^{\tip\wedge \theta}     |Y_s- Y_{\phi(s)}|^2ds } 
        \nonumber\\
	&+& 
	(C_L\;\alpha^{-1}+C^\alpha_L\;h)
	\Esp{
	\int_{t\wedge \theta}^{\tip\wedge \theta}\left(h+ \|Z_s-\hat Z_{\phi(s)}\|^2+\|\hat Z_{\phi(s)}-\bar Z_{\phi(s)}\|^2 \right)  ds 
	} 
	\nonumber 
\ee
\be
	&+&  (C_L\;\alpha^{-1}+C^\alpha_L\;h)
	\Esp{ \int_{t\wedge \theta}^{\tip\wedge \theta} (\xi_L\1_{\tau\wedge \taub\le s\le \tau\vee\taub}  + \1_{\taub \le s<\tau}   \|Z_s\|^2 ) ds }\;. 
\label{eq borne Delta 2}
	\ee  
Plugging \reff{eq borne Delta 2}  in \reff{eq borne Delta 1}  applied with $t=\ti$, using  Remark \ref{rem dble int},  taking $\alpha>0$ large enough, depending on the constants $C_L$, and $h$ small leads to 
	\b*
	\Delta^\theta_{\ti,\tip}
	&\le&
	(1+C_L\;h) \Esp{ |Y_{\tip\wedge\theta}-\bar Y_{\tip\wedge\theta}|^2 } 
	\nonumber\\
	&+& C_L\;\Esp{\int_{\ti\wedge \theta}^{\tip\wedge \theta}   \left( h+   |Y_s- Y_{\phi(s)}|^2+\|Z_s-\hat Z_{\phi(s)}\|^2 \right)  ds } 
	\nonumber\\
	&+&   C_L\; \Esp{ \int_{\ti\wedge \theta}^{\tip\wedge \theta} (\xi_L\1_{\tau\wedge \taub\le s\le \tau\vee\taub} + \1_{\taub \le s<\tau}   \|Z_s\|^2 ) ds }
	 \;.
	\e*
This implies that
	\b*
	\Delta^\theta
	&:=&
	\max_{i<n}\Esp{|Y_{\ti \wedge \theta}-\bar Y_{\ti \wedge \theta}|^2}
	+
	\Esp{\int_0^\theta \|Z_s-\tilde Z_s\|^2 ds}
	\\
	&\le&
	C_L\;\left(\Esp{ |Y_\theta-\bar Y_\theta|^2} 
	+h+ \Rc(Y)_{\Sc^2}^\pi+ \Rc(Z)_{\Hc^2}^\pi \right)
	\\
	&+& C_L\; \Esp{\xi_L\;|\taub\wedge \theta-\tau \wedge\theta|+ \int_{0}^{\theta} \1_{\taub \le s<\tau}   \|Z_s\|^2  ds } \;. 
	\e* 
We conclude the proof by using Remark \ref{rem dble int}  again to obtain
	\be
	\Esp{\int_0^\theta   \|Z_s-\bar Z_{\phi(s)}\|^2}
	&\le&
	C_L\;\left( \Esp{\int_0^\theta   \|\hat Z_{\phi(s)}-\bar Z_{\phi(s)}\|^2ds}
	+
	\Esp{\int_0^T   \|Z_s-\hat Z_{\phi(s)}\|^2ds}
	\right)
	\nonumber \\
	&\le&
	C_L\;\left(\Esp{\int_0^\theta   \|Z_s-\tilde Z_{s}\|^2ds}
	+
	\Esp{\int_0^T   \|Z_s-\hat Z_{\phi(s)}\|^2ds}\right) \label{eq cont Z bar Zphi par Z hat Z et Z tilde Z}
	\ee
which implies the required result, by the definition of Err$(h)^2_\theta$ in \reff{eq def Err}.\ep
\\

The above result  implies the first estimate of  Proposition \ref{prop majo erreur en terme de regu}.

\vs2

\no {\bf Proof of \reff{eq majo erreur en terme de regu T} of  Proposition \ref{prop majo erreur en terme de regu}.}  
It suffices to apply Proposition \ref{prop erreur BSDE theta} for $\theta=T$ and observe that the Lipschitz continuity of $g$  implies that 
	\b*
	\Esp{|g(\tau,X_\tau)-g(\taub,\Xb_\taub)|^2}
	&\le&
	C_L\;
	\Esp{
	|\tau-\taub|^2 
	+ 
	\|X_\taub-\Xb_\taub\|^2
	+
	\|\int_{\tau\wedge \taub}^{\tau\vee \taub} b(X_s) ds+\int_{\tau\wedge \taub}^{\tau\vee \taub} \sigma(X_s) dW_s\|^2 
	}
	\e*
where $|\tau-\taub|^2\le T|\tau-\taub|$, $\Esp{\|X_\taub-\Xb_\taub\|^2}\le	C_Lh$ by Proposition \ref{prop X bar X}, and 
	\b*
	\Esp{
	\|\int_{\tau\wedge \taub}^{\tau\vee \taub} b(X_s) ds+\int_{\tau\wedge \taub}^{\tau\vee \taub} \sigma(X_s) dW_s\|^2 
	}
	&\le&
	\Esp{\xi_L|\tau-\taub|}
	\e*
by Doob's inequality, \A{HL} and Proposition \ref{prop X bar X} again. 
\ep
\\

In order to prove \reff{eq majo erreur en terme de regu tau+ taub} of  Proposition \ref{prop majo erreur en terme de regu}, we  need the following easy Lemma. 

\begin{Lemma}\label{lem borne bar Y bar Z} Let \A{HL} hold. Then, 
	\be\label{eq majo bar Y bar Z tilde Z}
	\max_{i<n}\left( \|\bar Y_\ti\| + \sqrt{h} \|\bar Z_\ti\|\right) \le  \xi_L \And  \|\bar Y\|_{\Sc^2}+\|\bar Z_{\phi}\|_{\Hc^2}+\|\tilde Z\|_{\Hc^2} \le C_L\;.
	\ee
\end{Lemma}  
 
\proof The first bound follows from  the same arguments as in the proof of Lemma 3.3 in \cite{BoTo04}, after noticing that the boundedness assumption on $b$ and $\sigma$ can be relaxed for our result. Since, by \reff{eq def Euler BSDE avec tilde Z},
	\b*
	\bar Y_t=\Esp{\bar Y_\tip~|~\Fc_t}+\1_{\ti<\taub}(\tip-t)f(\Xb_\ti,\bar Y_\ti,\bar Z_\ti)
	\e*
on $[\ti,\tip]$, 
combining Jensen's inequality with \A{HL}, the first inequality of \reff{eq majo bar Y bar Z tilde Z} and Proposition \ref{prop X bar X}  imply that 	
	\be
	\sup_{t\le T} \Esp{|\bar Y_t|^2}
	&\le& 2\max_{i<n}\Esp{ |\bar Y_\tip|^2}
	+2
	h^{2} \max_{i\le n}\Esp{f(\Xb_\ti,\bar Y_\ti,\bar Z_\ti)^2}
	\sle C_L\;.\label{eq borne sup esp bar Y}
	\ee
Applying Itô's Lemma to $\bar Y^2$, using the inequality $ab\le a^2+b^2$ for $a,b\in \R$, \A{HL}, \reff{eq borne sup esp bar Y} and Proposition \ref{prop X bar X} then leads to
	\b*
	\Esp{\bar Y_{t\wedge \taub}^2}+\Esp{\int_{t\wedge \taub}^\taub \|\tilde Z_s\|^2ds}
	&=&
	\Esp{g(\taub,\Xb_\taub)^2+\int_{t\wedge \taub}^\taub 2\bar Y_s f(\bar X_{\phi(s)},\bar Y_{\phi(s)},\bar Z_{\phi(s)})  ds}
	\\
	&\le&
	C_L\;\left(1+ \alpha+\alpha^{-1}+ \alpha^{-1}\Esp{ \int_{t\wedge \taub}^\taub \|\bar Z_{\phi(s)}\|^2 ds }\right)\;,
	\e*
for all $\alpha>0$. By Remark \ref{rem dble int}, this shows that
	\b*
	\Esp{\int_{0}^\taub \|\bar Z_{\phi(s)}\|^2ds}\le\Esp{\int_{0}^\taub \|\tilde Z_s\|^2ds}
	&\le&
	C_L\;\left(1+ \alpha+\alpha^{-1}+ \alpha^{-1}\Esp{ \int_{0}^\taub \|\tilde Z_{s}\|^2 ds }\right)\;.
	\e*
Thus, taking $\alpha$ large enough, but depending only on $L$, and recalling Remark \ref{rem Z=0} leads to the required bound for  $\|\tilde Z\|_{\Hc^2}$ and $\|\bar Z_{\phi}\|_{\Hc^2}$. 
The bound on $\|\bar Y\|_{\Sc^2}$ is then easily deduced from its dynamics, Burkholder-Davis-Gundy's inequality, \A{HL}, \reff{eq borne sup esp bar Y} and Proposition \ref{prop X bar X}.\ep
\\

\no {\bf Proof of \reff{eq majo erreur en terme de regu tau+ taub} of  Proposition \ref{prop majo erreur en terme de regu}.}   Applying Proposition \ref{prop erreur BSDE theta} to $ \theta:=\tau_+ \wedge \taub$ and recalling Remark \ref{rem Z=0} leads to 
 	\b*
	{\rm {E}rr}(h)^2_{\tau_+\wedge \taub}
	&\le& 
	C_L\;\left(h+ \Esp{ |Y_{\tau_+\wedge \taub}-\bar Y_{\tau_+\wedge \taub}|^2} + \Rc(Y)_{\Sc^2}^\pi+ \Rc(Z)_{\Hc^2}^\pi \right) 
	 \;. 
	\e*
It remains to show that 
	\begin{equation}\label{eq borne bar Y Y tau wedge taub}
	\Esp{|\bar Y_{\tau_+\wedge \taub}-Y_{\tau_+\wedge \taub}|^2}  \le  C_L\left(h+ \Esp{\Esp{\xi_L|\tau-\taub|~|~\Fc_{\tau_+\wedge \taub}}^2}
	+\Esp{\1_{\taub<\tau}\Esp{ \int_\taub^\tau   \|Z_s\|   ds~|~\Fc_{\taub}}^2}\right)\;.
	\end{equation}

Since $f$ is $L$-Lipschitz continuous under \A{HL}, we can find an $\R^d$-valued adapted process $\chi$ which is bounded by $L$ and satisfies
	\be\label{eq def chi lipschitz}
	f(\bar X_{\phi(s)},\bar Y_{\phi(s)},\bar Z_{\phi(s)})=f(\bar X_{\phi(s)},\bar Y_{\phi(s)},0)+\scap{\chi_{\phi(s)}}{\bar Z_{\phi(s)}}
	\ee
on $[0,T]$. 
Set 
	\b*
	H_t:=\Ec\left(\int_0^t  \1_{\tau_+\le  s <  \taub} \chi_{\phi(s)} dW_s\right)\;,\;t\le T\;,
	\e*
where $\Ec$ stands for the usual Dol\'eans-Dade exponential martingale, 
and  define $\Q\sim \P$ by $d\Q/d\P=H_T$.  It follows from Girsanov's theorem that 
	\b*
	W^\Q=W - \int_0^\cdot \1_{\tau_+\le  s <  \taub} \chi_{\phi(s)}   ds 
	\e*
is a $\Q$-Brownian motion. Now, observe that, by \reff{eq def chi lipschitz} and \reff{eq def Euler BSDE avec tilde Z},
	\be
	Y_{t} &=& g(\tau,X_\tau)
	+ 
	\int_{t\wedge \tau}^\tau 
	f(X_{s},Y_s,Z_s) ds 
	-  
	\int_{t\wedge \tau}^\tau  Z_s dW^\Q_s\;
	 \label{eq def Y Q}
	\\
	\bar Y_{t} &=& g(\taub,\bar X_\taub)
	+ 
	\int_{t\wedge \taub}^\taub 
	\left(
	f(\bar X_{\phi(s)},\bar Y_{\phi(s)},\bar Z_{\phi(s)}) 
	-\1_{\tau_+\le  s}\scap{\chi_{\phi(s)}}{\tilde Z_{s}}
	\right)ds 
	-  
	\int_{t\wedge \taub}^\taub  \tilde Z_s dW^\Q_s\;. \label{eq def bar Y Q}
	\ee
In view of \reff{eq def chi lipschitz}, \reff{eq def Y Q}, \reff{eq def bar Y Q}, it then suffices to show that
	\be
	&&\Esp{\EspQ{g(\taub,\bar X_\taub) -g(\tau,X_{\tau})~|~\Fc_{\tau_+\wedge \taub}}^2} 
	\le 
	C_L\left(h+\Esp{\Esp{\xi_L|\tau-\taub|~|~\Fc_{\tau_+\wedge \taub} }^2}\right),\;~~~~~~~~~\label{eq prop control impact bsde sous Q eps control g}
	\\
	&&\Esp{\1_{\tau_+< \taub}\EspQ{
	\int_{\tau_+}^\taub f(\bar X_{\phi(s)},\bar Y_{\phi(s)},0)ds 
	~|~\Fc_{\tau_+}}^{2}}
	\le 
	\Esp{\Esp{\xi_L\;(|\tau-\taub|+h)~|~\Fc_{\tau_+\wedge \taub}}^{2}},~~~~~~~~~ 
	\label{eq prop control impact bsde sous Q eps 2}\\
	&&\Esp{\1_{\tau_+<\taub}\EspQ{
	\int_{\tau_+}^\taub \scap{\chi_{\phi(s)}}{\bar Z_{\phi(s)}-\tilde Z_{s}} ds 
	~|~\Fc_{\tau_+}}^{2}}\label{eq prop control impact bsde sous Q eps}
	\le C_L  h\;,
\ee
\be
	&&\Esp{\1_{\taub<\tau_+}\EspQ{
	\int_{\taub}^{\tau}  f(X_s,Y_s,Z_s) ds 
	~|~\Fc_{\taub}}^{2}}\nonumber
	\le C_L   \left(h+ \Esp{\Esp{\xi_L|\tau-\taub|~|~\Fc_{\tau_+\wedge \taub}}^2}\right)
	\\
	&&~~~~~~~~~~~~~~~~~~~~~~~~~~~~~~~~~~~~~~~~~~~~~~~~~~~~~~~~~~~ +C_L \Esp{\1_{\taub<\tau}\Esp{ \int_\taub^\tau   \|Z_s\|   ds~|~\Fc_{\taub}}^2}\;.\nonumber\\
	\label{eq prop control impact bsde sous Q eps 3}
	\ee
We start with the first term.   By using \A{HL}, applying   Itô's Lemma to $(g(t,X_t))_{t\ge 0}$ between $\taub$ and $\tau$, using Proposition \ref{prop X bar X}, the bound on $ \chi$ as well as standard estimates (recall \A{Hg} and Proposition \ref{prop X bar X}), we easily check that on $\{\tau_+> \taub\}\subset \{\tau>\taub\}$
	\b*
	\left|\EspQ{g(\tau, X_{\tau}) -g(\taub,\Xb_{\taub})~|~\Fc_{\taub}}\right|
	&\le&
	 C_L\;  \left\|X_{\taub}- \Xb_{\taub} \right\|
	 \\
	 &+&  \left|\EspQ{ \int_\taub^\tau  \left(\1_{\tau_+\le  s <  \taub} \scap{ \chi_{\phi(s)}\sigma^*}{Dg}   + \Lc g \right) (s,X_s) ds  ~|~\Fc_{\taub}} \right|
	 \\
	&\le& 
	 C_L\;   \left\|X_{\taub}- \Xb_{\taub} \right\|+   \Esp{\xi_L\;|\tau_+ -\taub|~|~\Fc_{\taub}} 
	\;.
	\e*
Similarly, on $\{\tau_+< \taub\}$, 
	\b*
	\left|\EspQ{g(\tau_+, X_{\tau_+}) -g(\taub,\Xb_{\taub})~|~\Fc_{\tau_+}}\right|
	&\le& 
	 C_L\;  \left\|X_{\tau_+}- \Xb_{\tau_+} \right\|+ \Esp{\xi_L\;|\tau_+-\taub|~|~\Fc_{\tau_+}}  
	\;.
	\e*
We then conclude the proof of \reff{eq prop control impact bsde sous Q eps control g} by appealing to \A{HL} and  Proposition \ref{prop X bar X} to obtain
	\b*
	\Esp{\left\|X_{\tau_+}- \Xb_{\tau_+} \right\|^2+ \left\|X_{\taub}- \Xb_{\taub} \right\|^2}+	\Esp{\left|g(\tau_+,X_{\tau_+})-g(\tau, X_\tau)\right|^2}
	 &\le& C_L  \; h \;,
	\e*
recall that $0\le \tau_+-\tau \le h$.	
	
The second term \reff{eq prop control impact bsde sous Q eps 2} is controlled by appealing to  \A{HL}, Lemma \ref{lem borne bar Y bar Z}  and Proposition \ref{prop X bar X}, recall that $\tau_+-\tau\le h$.
Concerning the third term \reff{eq prop control impact bsde sous Q eps}, we observe that $\{\tau_+\le s\}= \{ \tau \le \phi(s)\} \in \Fc_{\phi(s)}$ and that $\{\taub>s\}=\{\taub>\phi(s)\} \in \Fc_{\phi(s)}$. It then follows from \reff{eq lien bar Z tilde Z} that, on $\{\tau_+<\taub\}$,  
	\b*
	&&\EspQ{
	\int_{\tau_+}^\taub \scap{\chi_{\phi(s)}}{\bar Z_{\phi(s)}-\tilde Z_{s}} ds 
	~|~\Fc_{ \tau_+\wedge \taub}}
	\\
	&&= 
	\Esp{
	\int_{\tau_+}^\taub H_s \scap{\chi_{\phi(s)}}{\bar Z_{\phi(s)}-\tilde Z_{s}} ds 
	~|~\Fc_{\tau_+}}
	\\
	&&=
	\Esp{
	\int_{\tau_+}^\taub\; H_{\phi(s)}
	\left\langle \chi_{\phi(s)}\;,\;h^{-1}\int_{\phi(s)}^{\phi(s)+h} \tilde Z_{u}du-\tilde Z_{s}\right\rangle ds 
	~|~\Fc_{\tau_+}}
	\\
	&&+
	\Esp{
	\int_{\tau_+}^\taub \; (H_s-H_{\phi(s)})
	\left\langle \chi_{\phi(s)}\;, \bar Z_{\phi(s)}-\tilde Z_{s}\right\rangle ds 
	~|~\Fc_{\tau_+}}
	\e*
and, since $\taub$ and $\tau_+$ take values in $\pi$, 
	\b*
	\int_{\tau_+}^\taub\; H_{\phi(s)}
	\left\langle \chi_{\phi(s)}\;,\;h^{-1}\int_{\phi(s)}^{\phi(s)+h} \tilde Z_{u}du-\tilde Z_{s}\right\rangle ds 
	&=&0\;. 
	\e*
On the other hand, the Cauchy-Schwartz inequality and the boundedness of $\chi$ imply that 
	\b*
	&&
	\left|\Esp{
	\int_{\tau_+}^\taub \; (H_s-H_{\phi(s)})
	\left\langle \chi_{\phi(s)}\;,\;\bar Z_{\phi(s)}-\tilde Z_{s}\right\rangle ds 
	~|~\Fc_{\tau_+\wedge \taub}}\right|
	\\
	&&
	\sle C_L\;
	\left|\Esp{
	\int_{\tau_+}^\taub \; (H_s-H_{\phi(s)})^2 ds 
	~|~\Fc_{\tau_+\wedge \taub}}\right|^\frac12\;
	\left|\Esp{
	\int_{\tau_+}^\taub \; \|\bar Z_{\phi(s)}-\tilde Z_s\|^2 ds 
	~|~\Fc_{\tau_+\wedge \taub}}\right|^\frac12
	\\
	&&
	\le
	\xi_L h^{\frac12}\left|\Esp{
	\int_{\tau_+}^\taub \; \|\bar Z_{\phi(s)}-\tilde Z_s\|^2 ds 
	~|~\Fc_{\tau_+\wedge \taub}}\right|^\frac12\;. 
	\e*
Recalling Lemma \ref{lem borne bar Y bar Z} and combining the above inequalities leads to  \reff{eq prop control impact bsde sous Q eps}. 

The last term \reff{eq prop control impact bsde sous Q eps 3} is easily controlled by using \A{HL}, Remark \ref{rem Z=0},  and Proposition \ref{prop X bar X}. \ep
 

\section{Exit time approximation error: Proof of Theorem \ref{thm main erreur ta}}\label{sec thm main erreur ta}

In this section, we provide the proof of  Theorem \ref{thm main erreur ta}.  We start with a partial argument which essentially allows to reduce to the case where $m=1$, i.e. $\Oc$ has no corners, by working separately on the exit times of the different domains $\Oc^\ell$:  
	\b*
	\tau^\ell_+ := \inf\{t\in \pi~:~ \exists\;s\le t ~{\rm s.t. }~ X_{s}\notin \Oc^\ell\}\wedge T\;
	& \mbox{ and }& 
	\taub^\ell := \inf\{t\in \pi~:~\bar X_{t}\notin \Oc^\ell\}\wedge T\;\;.
	\e*
We shall prove below the following Proposition.

\begin{Proposition}\label{prop thm main erreur ta ell} Assume that \A{HL}, \A{D1} and \A{C} hold. Then, for each $\eps>0$, 
	\be\label{eq thm main erreur ta ell} 
	\Esp{ \Esp{|\tau^\ell_+ -\taub^\ell|~|~\Fc_{\tau^\ell_+\wedge \taub^\ell} }^2} &\le& C_L^\eps h^{1-\eps},\; \;\;\;\forall\;1\le \ell \le m\;. 
	\ee
\end{Proposition} 

It  implies the statements of  Theorem \ref{thm main erreur ta}.

\vs2

\no {\bf Proof of Theorem \ref{thm main erreur ta}.} 
Since  $\tau_+=\min_{\ell \le m} \tau^\ell_+$  and $\taub=\min_{\ell \le m} \taub^\ell$, 
we have
	\b*
	\Esp{|\tau_+ -\taub|~|~\Fc_{\tau_+\wedge \taub} }
	&\le&
	\sum_{\ell=1}^m
	\Esp{|\tau^\ell_+ -\taub^\ell|~|~\Fc_{\tau^\ell_+\wedge \taub^\ell} }\left(\1_{\tau_+=\tau_+^\ell<\taub}+\1_{\taub=\taub^\ell\le\tau_+}\right)
	\e*
which combined with   \reff{eq thm main erreur ta ell} leads to 
	\be\label{eq proof thm main erreur ta}
	\Esp{ \Esp{|\tau -\taub|~|~\Fc_{\tau_+\wedge \taub} }^2} &\le& C_L^\eps h^{1-\eps}\;, 
	\ee
since $|\tau_+ -\tau|\le h$. This leads to the second assertion of  Theorem \ref{thm main erreur ta}. 
On the other hand, given a positive random variable $\xi$ satisfying $\Esp{\xi^p}\le C_L^p$ for all $p\ge 1$, we deduce from Hölder's inequality that 
   	\b*
   	\Esp{\xi \;|\tau -\taub|~|~\Fc_{\tau_+\wedge \taub} }^2
   	&\le&
   	\xi_L^\eps\; \Esp{|\tau -\taub|^\frac{1}{1-\eps}~|~\Fc_{\tau_+\wedge \taub} }^{2(1-\eps)}
   	\sle
   	\xi_L^\eps\; T^{2\eps}\;\Esp{|\tau -\taub|~|~\Fc_{\tau_+\wedge \taub} }^{2(1-\eps)}
   	\e*
and 
		\b*
   	\Esp{\xi\;\Esp{\xi \;|\tau -\taub|~|~\Fc_{\tau_+\wedge \taub} }^2}
   	&\le&
   	C_L^\eps\; \Esp{\Esp{|\tau -\taub|~|~\Fc_{\tau_+\wedge \taub} }^{2}}^{1-\eps}\;.
   	\e*
In view of \reff{eq proof thm main erreur ta}, this leads to the first assertion of  Theorem \ref{thm main erreur ta}, after possibly changing $\eps$.  \ep

\vs2

The rest of this section is devoted to the proof of \reff{eq thm main erreur ta ell} for some fixed $\ell$. 
We first provide an \emph{a-priori} control on the difference between $\tau^\ell_+$ and $\taub^\ell$. We use the standard idea that consists in introducing a test function on which we can apply Itô's Lemma between  $\tau^\ell_+$ and $\taub^\ell$ so that the Lebesgue integral term provides an upper bound for the difference between these two times, see e.g. Lemma 3.1 Chapter 3 in \cite{Fr85} for an application to the construction of upper bounds for the moments of the first exit time of a  uniformly elliptic diffusion from a bounded domain. 

\vs2

To this end, we introduce the family of test functions
	$$
	F_\ell:=d_\ell^2/\gamma\;\;,\;1\le \ell \le m\;,
	$$ 
for some $\gamma>0$ to be fixed below. Here,  $d_\ell$ is a $C^2(\R^d)$ function which coincides with the algebraic distance to $\partial \Oc^\ell$ on a neighborhood of $\partial \Oc^\ell$ and such that 
	\b*
	\Oc^\ell:=\{x\in \R^d~:~d_\ell(x)>0\} &\And& \partial \Oc^\ell:=\{x\in \R^d~:~d_\ell(x)=0\}\;. 
	\e* 
The existence of such a map is guaranteed by the smoothness assumption \A{D1}, see e.g. \cite{GiTr98}. Observe that, after possibly changing $L$ and considering a suitable extension of $d_\ell$ outside of a neighbourhood of the compact boundary $\partial \Oc^\ell$, we can assume that 
	\be\label{eq borne derivees d}
	\|d_\ell\|+\|Dd_\ell\|+\|D^2d_\ell\|&\le& L\;\; \mbox{ on } \R^d\;. 
	\ee 	
Observe that
	\be
	\Lc F_\ell
	&=&
	\frac{1}{\gamma} \left[\left(2 \scap{b}{n_\ell}+ \Tr{ a D^2 d_\ell }\right)d_\ell + \Tr{ a (n_\ell)^* n_\ell }\right]
	\ee
where 
$
	n_\ell:=D d_\ell\;
	$
coincides with the unit inward normal for $x\in \partial \Oc^\ell$, recall \A{D1}. 

In view of \A{HL}, \A{D1}, \reff{eq borne derivees d} and \A{C}, there is some $C_L>0$ such that, for each  $1\le \ell \le m$, 
\be\label{eq LcF ge 1}
	\Lc F_\ell\ge \frac{1}{\gamma} ( -C_L d_\ell +  n_\ell   \;a (n_\ell)^* )
	\ge 1 \mbox{ and } n_\ell  \; a (n_\ell)^*  \ge L^{-1}/2
  	&\mbox{ on }& B(\partial \Oc^\ell,r)\;\;
	\ee
if we choose $r>0$ and $\gamma>0$ small enough, but depending only on $L$. For later use, also observe that, after possibly changing $r$, one can actually choose it such  that 
	\be\label{eq LcF ge 1 bis}
	 n_\ell(x) \;a(y) n_\ell(x)^*  \ge L^{-1}/2
  	&\mbox{ for all }& x,y  \in B(\partial \Oc^\ell,r)\;\;{\rm s.t. } \;\|x-y\|\le r\;\;. 
	\ee

\vs2
 
We now fix $r,\gamma>0$ such that \reff{eq LcF ge 1} and \reff{eq LcF ge 1 bis} hold and define the sets
	\b*
	A_\ell&:=&\{X_s \in B(\partial \Oc^\ell,r)\;,\;\forall\;s\in [\bar \tau^\ell,\tau^\ell_+]\} \;,\; B_\ell:= \{|d_\ell(X_{\tau^\ell_+})|\le h^{\frac12-\eta}\} 
	\\
	\bar A_\ell&:=& \{\bar X_s \in B(\partial \Oc^\ell,r)\;,\;\forall\;s\in [\tau^\ell_+, \taub^\ell]\}
	\;,\; 
	\bar B_\ell:=\{|d_\ell(\bar X_{\taub^\ell})|\le h^{\frac12-\eta}\} \;
 	\;,
	\e*
for some $\eta \in (0,1/4)$ to be chosen later on. Observe that $A_\ell$ (resp. $\bar A_\ell$) is well defined   on $\{\bar \tau^\ell\le \tau^\ell_+\}$ (resp.  $\{\tau^\ell_+\le \taub^\ell\}$).  

\vs2

We can now provide our first control on $|\tau^\ell_+- \taub^\ell|$. 
Recall that  $\xi_L^\eps $ ($\xi_L$ if it does not depend on some extra parameter $\eps$) denotes a positive random variable whose value may change from line to line but  satisfies $\Esp{|\xi_L^\eps|^p}\le C_L^{\eps,p}$ for all $p\ge 1$.

\begin{Lemma}\label{lem borne ecart tau borne par PA} Assume that \A{HL} and \A{D1} hold. Then, for each $\eps\in (0,1)$, 
	\b*
	\Esp{ |\tau^\ell_+- \taub^\ell| ~|~\Fc_ {\tau^\ell_+ \wedge \taub^\ell} }
	&\le& 
	\xi_L^\eps\left\{ h^{\frac12} 
	+ \;(T-\taub^\ell)^\frac12 \Pro{ (A_\ell \cap B_\ell)^c ~|~\Fc_ {\taub^\ell} }^{1-\eps} \1_{\{\tau^\ell_+> \taub^\ell\}}   
	\right. \\
	&&
	\;\;\;\;\;\;\;\;\;\;\;\;+ \;
	\left.
	(T-\tau^\ell_+)^\frac12
	\Pro{ (\bar A_\ell \cap \bar B_\ell)^c ~|~\Fc_ {\tau^\ell_+ } }^{1-\eps} \1_{\{\tau^\ell_+< \taub^\ell\}} 
	  \right\}\; 
	\e*
for each $1\le \ell\le m$.
\end{Lemma}

\proof 1. We first work on the event  $\{\tau^\ell_+>\taub^\ell\}$.  It   follows from  \reff{eq LcF ge 1}  and Itô's Lemma that  
	\b*
	\Esp{ \tau^\ell_+- \taub^\ell ~|~\Fc_ {\taub^\ell}}
	&\le&
	\Esp{
	\1_{A_\ell\cap B_\ell} \int_{\taub^\ell}^{\tau^\ell_+} \Lc F_\ell(X_s) ds ~|~\Fc_ {\taub^\ell}} + (T-\taub^\ell)\; \Pro{ (A_\ell\cap B_\ell)^c  
	~|~\Fc_ {\taub^\ell} }
	\\
	&\le&
	\Esp{
	\1_{A_\ell\cap B_\ell}  \left(\int_{\taub^\ell}^{\tau^\ell_+} \Lc F_\ell(X_s) ds  + \int_{\taub^\ell}^{\tau^\ell_+} D F_\ell(X_s)\sigma(X_s) dW_s\right)
	~|~\Fc_ {\taub^\ell}} \\
	&-&
	\Esp{
	\1_{A_\ell\cap B_\ell}  \int_{\taub^\ell}^{\tau^\ell_+} D F_\ell(X_s)\sigma(X_s) dW_s 
	~|~\Fc_ {\taub^\ell}}
	\\
	&+& (T-\taub^\ell)\; \Pro{ (A_\ell\cap B_\ell)^c  ~|~\Fc_ {\taub^\ell} }
	\\
	&\le&
	\gamma^{-1} \;\Esp{ (d_{\ell}^2( X_{\tau^\ell_+})-d_{\ell}^2(X_{\taub^\ell}))  \1_{A_\ell\cap B_\ell}~|~\Fc_ {\taub^\ell}}
	\\
	&+& 
	\Esp{\1_{(A_\ell\cap B_\ell)^c }  \int_{\taub^\ell}^{\tau^\ell_+} D F_\ell(X_s)\sigma(X_s) dW_s ~|~\Fc_ {\taub^\ell} }
	\\
	&+& (T-\taub^\ell)\; \Pro{ (A_\ell\cap B_\ell)^c  ~|~\Fc_ {\taub^\ell} }
	\e*
where, by  Hölder's and Burkholder-Davis-Gundy's inequality, the Lipschitz continuity of $\sigma$ and $DF_\ell$ (see \A{HL} and \reff{eq borne derivees d}) and Proposition \ref{prop X bar X}, 
	\b*
	\Esp{\1_{(A_\ell\cap B_\ell)^c }  \int_{\taub^\ell}^{\tau^\ell_+} D F_\ell(X_s)\sigma(X_s) dW_s ~|~\Fc_ {\taub^\ell} }
	&\le& 
	\xi_L^{\eps}\;(T-\taub^\ell)^\frac12 \Pro{(A_\ell\cap B_\ell)^c ~|~\Fc_ {\taub^\ell} }^{1-\eps}  \; 
	\e*
for all $\eps \in (0,1)$. 
We now recall that  $|d_\ell(X_{\tau^\ell_+})|$ $\le$ $h^{\frac12-\eta}$  on $B_\ell$, which implies 
	\b*
	\Esp{( d_{\ell}^2(X_{\tau^\ell_+})-d_{\ell}^2(X_{\taub^\ell})) \1_{A_\ell\cap B_\ell}~|~\Fc_ { \taub^\ell} } \le \Esp{  d_\ell^2(X_{\tau^\ell_+})  \1_{A_\ell\cap B_\ell}~|~\Fc_ { \taub^\ell} } 
	&\le&    h^{1-2\eta} \;.
	\e*
In view of the above inequalities, this provides the required estimate  on the event set  $\{\tau^\ell_+>\taub^\ell\}$ since   $\eta<1/4$.
	  
2. We now work on  the event    $\{\tau^\ell_+<\taub^\ell\}$.   By Proposition \ref{prop X bar X},
	\b*
	\Esp{\1_{\bar A_\ell \cap \bar B_\ell}  
	\int_{\tau^\ell_+}^{\taub^\ell} \left| \Lc^{\bar X_{\phi(s)}} F_\ell(\bar X_s)-\Lc^{\bar X_{s}} F_\ell(\bar X_s)\right| ds 
	~|~\Fc_{\tau^\ell_+}}
	\le \xi_L \;h^{\frac12}\;,
	\e*
with the notation
$	\Lc^y F_\ell := \partial_t F_\ell+  \scap{b(y)}{D F_\ell}  +  \frac12  \Tr{ a(y) D^2 F_\ell   } \;$, 
so that  $\Lc^{\bar X_{s}} F_\ell(\bar X_s)=\Lc F_\ell(\bar X_s)$.
Arguing as above, it follows that, on  $\{\taub^\ell>\tau^\ell_+\}$, 
	\b*
	\Esp{ \taub^\ell- \tau^\ell_+ ~|~\Fc_{\tau^\ell_+}}
	&\le&
	\xi_L \;h^{\frac12}
	+
	\gamma^{-1} \;\Esp{ (d^2_\ell( \bar X_{\taub^\ell})-d^2_\ell(\bar X_{\tau^\ell_+}))  \1_{\bar A_\ell\cap \bar B_\ell}~|~\Fc_ {\tau^\ell_+}}
	\\
	&+& 
	\Esp{\1_{(\bar A_\ell\cap \bar B_\ell)^c }  \int_{\tau^\ell_+}^{\taub^\ell} D F_\ell(\bar X_s)\sigma(\bar X_{\phi(s)}) dW_s ~|~\Fc_ {\tau^\ell_+} }
	\\
	&+& (T-\tau^\ell_+)\; \Pro{ (\bar A_\ell\cap \bar B_\ell)^c  ~|~\Fc_ {\tau^\ell_+} }
	\\
	&\le&\xi_L \;h^{\frac12}
	+
	\gamma^{-1} \;h^{\frac12}  
	+\xi_L^\eps\;(T-\tau^\ell_+)^\frac12 \Pro{(\bar A_\ell\cap \bar B_\ell)^c ~|~\Fc_ {\tau^\ell_+} }^{1-\eps} 
	\;. 
	\e*
\ep
\\

It remains to control the different terms that appear in the upper bound of Lemma \ref{lem borne ecart tau borne par PA}.

\vs2

For notational convenience, we now introduce the sets (recall that $0<\eta<1/4$)
	\b*
	E_\ell:=\{d_\ell(X_{\taub^\ell})  \le  h^{\frac12 -\eta}  \} 
	&\And& \bar E_\ell  := \{d_\ell(\bar X_{\tau^\ell_+})  \le  h^{\frac12 -\eta}\}\;\;,\;1\le \ell \le m\;.
	\e*

\begin{Remark}\label{rem Pro Ec}{\rm Observe that 
	\b*
	\Pro{E_\ell^c\cap \{ \bar \tau^\ell<\tau^\ell_+\}}
	&\le&
	\Pro{E_\ell^c \cap \{\bar \tau^\ell<T\}}
	\sle
	\Pro{\{d_\ell(X_{\taub^\ell})-d_\ell(\Xb_{\taub^\ell})\ge  h^{\frac12 -\eta}\} \cap \{\bar \tau^\ell<T\}}\;,
	\e*
since $d_\ell(\Xb_{\taub^\ell})\le 0$ on $\{\bar \tau^\ell<T\}$. 
Using \reff{eq borne derivees d}, Tchebychev's inequality and Proposition \ref{prop X bar X}, we then deduce that, for each $\eps \in (0,1)$, there is $C_L^\eps>0$ such that 
	\b*
	\Pro{E_\ell^c\cap \{ \bar \tau^\ell<\tau^\ell_+\}}
	&\le&
	C^\eps_L\;h^{1-\eps}\;.
	\e*

Similarly, if $\tau^\ell$ denotes the first exit time of $(t,X_t)_{t\ge 0}$ from $[0,T)\x \Oc^\ell$, we have
	\b*
	\Pro{\bar E_\ell^c\cap \{ \bar \tau^\ell>\tau^\ell_+\}}
	&\le&
	\Pro{
	\{d_\ell(\Xb_{\tau^\ell_+})-d_\ell(X_{\tau^\ell_+})\ge  \frac12 h^{\frac12 -\eta}\}
	\cap 
	\{d_\ell(X_{\tau^\ell_+})\le  \frac12 h^{\frac12 -\eta}\} \cap \{\tau^\ell_+<T\}
	}
	\\
	&+& 
	\Pro{ \{d_\ell(X_{\tau^\ell_+})-d_\ell(X_{\tau^\ell})>   \frac12 h^{\frac12 -\eta}\} \cap \{\tau^\ell_+<T\}}
	\\
	&\le&
	C^\eps_L\;h^{1-\eps}\;,
	\e*
where the last inequality follows from Tchebychev's inequality, Proposition \ref{prop X bar X} and the fact that $\tau^\ell_+-\tau^\ell\le h$. Note that the term 
$d_\ell(X_{\tau^\ell_+})-d_\ell(X_{\tau^\ell})$ could be controlled by Bernstein type inequalities in order to avoid the explosion of the constant with $\eps$. However, to the best of our knowledge, such inequalities are not available in the existing literature for the term $d_\ell(\Xb_{\tau^\ell_+})-d_\ell(X_{\tau^\ell_+})$ and Tchebychev's inequality remains the most natural tool to apply here. 
}
\end{Remark}


Combining the above Remark with  the next two technical Lemmas allows to control  the right hand-side terms in the upper bound of Lemma \ref{lem borne ecart tau borne par PA}. Thus, the statement of Proposition \ref{prop thm main erreur ta ell} is a direct consequence of   Lemma \ref{lem borne ecart tau borne par PA} combined with Remark \ref{rem Pro Ec}, Lemma \ref{lem borne P[Ac] d(x)le sqrt(h)} and Lemma \ref{lem borne P[Bc cap A]} below, applied for $\eta$ small enough. 
	
\begin{Lemma}\label{lem borne P[Ac] d(x)le sqrt(h)} Assume that \A{HL}, \A{D1} and \A{C} hold. Then, for each $\eps \in (0,1)$,
	\be\label{eq lem borne P[Ac] d(x)le sqrt(h)}
	\Pro{A_\ell^c ~|~\Fc_{\taub^\ell}} \1_{E_\ell\cap \{\tau^\ell_+>\bar \tau^\ell\}}
	+ \Pro{\bar A_\ell^c  ~|~\Fc_{\tau^\ell_+}}\1_{\bar E_\ell\cap \{\tau^\ell_+< \taub^\ell\}}
	&\le& 
	\xi_L^\eps \; h^{(\frac12 -\eta)(1 -\eps)} \;,\;\forall \; \ell \le m \;. 
	\ee 
\end{Lemma}

\begin{Lemma}\label{lem borne P[Bc cap A]} Assume that \A{HL}, \A{D1} and \A{C} hold. Then, for each $\eps\in (0,1)$, 
		\begin{equation}\label{eq lem borne P[Bc cap A]}
		\Pro{A_\ell \cap B_\ell^c ~|~ \Fc_{\taub^\ell} }\1_{E_\ell\cap \{\tau^\ell_+>\taub^\ell\}} 
		+ 
		\Pro{\bar A_\ell \cap \bar B_\ell^c ~|~ \Fc_{\tau^\ell_+}} \1_{\bar E_\ell\cap \{\tau^\ell_+<\taub^\ell\}}
		  \le   
		\xi_L^\eps \; \frac{h^{(\frac12 -\eta)(1 -\eps)}}{\sqrt{T-\taub^\ell\wedge \tau^\ell_+}} \;,\;\forall \; \ell \le m \;. 
		\end{equation}
\end{Lemma}

\vs2
 
\no {\bf Proof of Lemma \ref{lem borne P[Ac] d(x)le sqrt(h)}.} 1. We first prove the bound for the first term.  
Let  $V$ be defined by $V_t:=d_\ell(X_{\taub^\ell+t})$ for $t\ge 0$ and let $\vartheta^y$ be the first time when $V$ reaches $y \in \R$. Using $A_\ell^c=A_\ell^c\cap(\{\vartheta^0\ge \vartheta^r\}\cup \{\vartheta^0< \vartheta^r\})$, we deduce that on $\{\tau^\ell_+ >\taub^\ell\}\cap E_\ell$
	\b*
	\Pro{A_\ell^c  ~|~ \Fc_{\taub^\ell}}
	&\le& 
	\Pro{\vartheta^0\ge \vartheta^r ~|~ \Fc_{\taub^\ell}}
	+
	\Pro{ \{\sup_{s\in [\tau^\ell,\tau^\ell_+]} |d_\ell(X_s)|\ge r\}\cap\{\tau^\ell<T\} ~|~ \Fc_{\taub^\ell}}\;,
	\e*
where, by \reff{eq borne derivees d}, Tchebychev's inequality and Proposition \ref{prop X bar X}, on $\{\tau^\ell_+ >\taub^\ell\}\subset \{\tau^\ell>\taub^\ell\}$, 
	\b*
	\Pro{ \{\sup_{s\in [\tau^\ell,\tau^\ell_+]} |d_\ell(X_s)|\ge r\}\cap\{\tau^\ell<T\} ~|~ \Fc_{\taub^\ell}} 
	&\le&
	r^{-2} \Esp{\sup_{s\in [\tau^\ell,\tau^\ell_+]} |d_\ell(X_s)-d_\ell(X_{\tau^\ell})|^2~|~ \Fc_{\taub^\ell}} 
	\\
	&\le&
	\xi_L\;h\;,
	\e*
recall that $\tau^\ell_+-\tau^\ell\le h$.
It remains to provide a suitable bound for $\Pro{\vartheta^0\ge \vartheta^r ~|~ \Fc_{\taub^\ell}}$. From now on, we assume, without loss of generality, that 
	\be\label{eq 2hler}
	2h^{\frac12-\eta} \le r\;.
	\ee

Set $\vartheta:=\vartheta^0 \wedge \vartheta^r$. Thanks to \A{C} and \A{HL}, we can define   $\Q\sim \P$  by the density 
	\b*
	H=\Ec_{\taub^\ell+\vartheta} \left( 
	-\int_0^\cdot \1_{E_\ell}\1_{s\ge {\taub^\ell}}  (n_\ell\sigma) (X_s)((n_\ell a n_\ell^*)(X_s))^{-1} \Lc d_\ell(X_s) dW_s
	\right) \;.
	\e*
Let 
	$$
	W^\Q:=W+\1_{[\taub^\ell,\infty)} \1_{E_\ell}
	\int_{\taub^\ell}^{(\taub^\ell+\vartheta)\wedge \cdot}   (n_\ell\sigma)^* (X_s)((n_\ell a n_\ell^*)(X_s))^{-1} \Lc d_\ell(X_s)ds
	$$ 
be the  Brownian motion associated to $\Q$ by Girsanov's Theorem. We have 
	\b*
	V_{t\wedge \vartheta}=V_0+ \int_{\taub^\ell}^{\taub^\ell+{t\wedge \vartheta}} n_\ell(X_s) \sigma(X_s) dW^\Q_s\;\;\;\mbox{ on } \; E_\ell\;. 
	\e*
Set 
	\b*
	\Lambda_t:= \int_{\taub^\ell}^{\taub^\ell+t} \|n_\ell(X_{s\wedge (\taub^\ell+\vartheta)}) \sigma(X_{s\wedge (\taub^\ell+\vartheta)})\|^2 ds\;. 
	\e*
By the Dambis-Dubins-Schwarz theorem, see Theorem 4.6 Chapter 3 in \cite{KaSh90}, there exists a one dimensional $\Q$-Brownian motion $Z$ such that 
	\b*
	V_{t\wedge \vartheta}=V_0+Z_{\Lambda_{t\wedge \vartheta}}\;\;\;\;
	\mbox{ on } \;E_\ell \cap \{\tau^\ell_+ >\taub^\ell\}=\{V_0\le   h^{\frac12 -\eta}\;,\;\tau^\ell_+ > \taub^\ell\}\;.
	\e*
This implies that
	\b*
	 \Qro{ \vartheta^0\ge \vartheta^r~|~ \Fc_{\taub^\ell} }
	 &\le&
	  h^{\frac12 -\eta} /r  \;\;\;\;\mbox{ on } \;E_\ell \cap \{\tau^\ell_+ >\taub^\ell\}\;, 
	\e*
see e.g. Exercise 8.13 Chapter 2.8 in \cite{KaSh90}.
We conclude by using Hölder's inequality and \reff{eq borne derivees d}. 

\no 2. The bound for the second term in \reff{eq lem borne P[Ac] d(x)le sqrt(h)} is derived similarly. We now write 
	\b*
	V_t:=d_\ell(\Xb_{\tau^\ell_+ + t})  \;,\;t\ge 0\;. 
	\e*
As above, we denote by $\vartheta^y$  the first time when $V$ reaches $y \in \R$ and observe that, by \reff{eq 2hler},
	\b*
	\Pro{\bar A_\ell^c  ~|~ \Fc_{\tau^\ell_+}}
	&\le&
	\Pro{\vartheta^{-h^{\frac12-\eta}}>\vartheta^r  ~|~ \Fc_{\tau^\ell_+}}
	+
	\Pro{\sup_{s\in [\tilde \tau^\ell,\tilde \tau^\ell+h]} |d_\ell(\Xb_{s})-d_\ell(\Xb_{\tilde \tau^\ell})|>h^{\frac12-\eta}  ~|~ \Fc_{\tau^\ell_+}}
	\e*
where $\tilde \tau^\ell:=\tau^\ell_+ +\vartheta^{-h^{\frac12-\eta}}$, and, by \reff{eq borne derivees d}, Tchebychev's inequality and Proposition \ref{prop X bar X},
	\b*
		\Pro{\sup_{s\in [\tilde \tau^\ell,\tilde \tau^\ell+h]} |d_\ell(\Xb_{s})-d_\ell(\Xb_{\tilde \tau^\ell})|>h^{\frac12-\eta}  ~|~ \Fc_{\tau^\ell_+}}
		&\le&
		\xi^\eta_L\;h\;.
	\e*
In order to bound the term $\Pro{\vartheta^{-h^{\frac12-\eta}}>\vartheta^r  ~|~ \Fc_{\tau^\ell_+}}$, we observe that \reff{eq LcF ge 1 bis}   imply that, for $h$ small enough, 
	\b*
	\|n_\ell(\bar X_s)\sigma(\bar X_{\phi(s)})\|\ge  L^{-\frac12}/\sqrt{2} 
	\;\;\mbox{ on } \;\bar E_\ell\cap\{s \in [\tau^\ell_+,\theta^\ell]\}\cap \{ \|\bar X_s-\bar X_{\phi(s)}\|\le r\}\;,
	\e*
where $\theta^\ell:=\inf\{t\ge \tau^\ell_+~:~\Xb_t \notin B(\partial \Oc^\ell,r)\}\wedge T$.  Moreover, it follows from Proposition \ref{prop X bar X} that 
	\b*
	\Pro{\sup_{s\le T} \|\bar X_s-\bar X_{\phi(s)}\|> r} &\le& C_L\;r^{-4}\; h \;.
	\e*
Up to obvious modifications, this allows us to reproduce the arguments of Step 1 on  the event set $\bar E_\ell$.
\ep  
\\

\no {\bf Proof of Lemma \ref{lem borne P[Bc cap A]}.} We only prove the bound for the first term. The second one can be derived from similar arguments (see step 2 in the proof of Lemma \ref{lem borne P[Ac] d(x)le sqrt(h)}). We use the notations of the proof of Lemma \ref{lem borne P[Ac] d(x)le sqrt(h)}. We first observe that, on $E_l\cap \{\tau^\ell>\bar \tau^\ell \} $,
	\b*
	\Pro{A_\ell \cap B_\ell^c ~|~ \Fc_{\taub^\ell}}  
	&\le& 
	\Pro{A_\ell\cap\{\vartheta^0>(T-\taub^\ell)\}~|~\Fc_{\taub^\ell}}\;
	\\
	&+&
	\Pro{\{\tau^\ell<T\}\cap \sup_{s\in [\tau^\ell,\tau^\ell_+]} |d_\ell(X_s)-d_\ell(X_{\tau^\ell})|\ge h^{\frac12-\eta}|~|~\Fc_{\taub^\ell}}\;
	\\
	&\le&
  \Pro{A_\ell\cap\{\min_{t\in [0,T-\taub^\ell]}  Z_{\Lambda_t}  >-h^{\frac12-\eta}\} ~|~\Fc_{\taub^\ell}}+\xi^\eta_L\;h\;
   \;,
  \e*
where the second inequality follows from Tchebychev's inequality, \A{HL} and Proposition \ref{prop X bar X}, recall that $\tau^\ell_+-\tau^\ell \le h$.
Using Hölder's inequality, we then observe that 
  \b*
  \Pro{A_\ell\cap\{\min_{t\in [0,T-\taub^\ell]}  Z_{\Lambda_t}  >-h^{\frac12-\eta}\} ~|~\Fc_{\taub^\ell}}\;
  &\le&   
  \xi_L^\eps\; \Qro{A_\ell\cap\{\min_{t\in [0,T-\taub^\ell]} Z_{\Lambda_t} >-h^{\frac12 -\eta}\}~|~\Fc_{\taub^\ell}}^{1-\eps}\;.  
	\e*
Since, by \reff{eq LcF ge 1 bis}, 
	\b*
	\Lambda_{T-\taub^\ell}\ge  (T-\taub^\ell)(2L)^{-1}\;\;\;
	\mbox{on}\;
	A_\ell \cap\{\vartheta^0>(T-\taub^\ell)\}\cap\{\taub^\ell< \tau^\ell_+\}\subset 	A_\ell \cap\{\taub^\ell< \tau^\ell_+=T \}\;,
	\e* 
we deduce from Chapter 2 of \cite{KaSh90} that, on $E_\ell\cap\{\taub^\ell< \tau^\ell_+\}$, 
		\b*
	  \Qro{A_\ell\cap\{\min_{t\in [0,T-\taub^\ell]} Z_{\Lambda_t} >-h^{\frac12 -\eta}\}~|~\Fc_{\taub^\ell}} 		 
	  &\le&
	  \Qro{\min_{t\in [0,(T-\taub^\ell)(2L)^{-1}]}  Z_{t}  > -h^{\frac12 -\eta}~|~\Fc_{\taub^\ell}}	
	  \\
	  &\le&
	  C_L\;(T-\taub^\ell)^{-\frac12} h^{\frac12-\eta}\;. 
	\e*
We conclude by combining the above estimates. 
\ep

\section{Regularity of the BSDE and the related PDE}\label{sec regu}

 \subsection{Interpretation in terms of parabolic semilinear PDEs with Dirichlet boundary conditions}\label{sec proof of Proposition prop u PDE}

In this section, we denote by $X^{t,x}$ the solution of \reff{eq def SDE intro} with initial condition $x \in \bar \Oc$ at time $t\le T$. We also denote by $\tau^{t,x}$ the first exit time of $(s,X^{t,x}_s)_{s\ge t}$ from $\Oc\x [0,T)$ and write $(Y^{t,x},Z^{t,x})$ for the solution of \reff{eq def BSDE intro} with $(X^{t,x},\tau^{t,x})$ in place of $(X,\tau)$.

As usual the deterministic function $(t,x)\in \bar D \mapsto u(t,x):=Y^{t,x}_t$ can be related to the semilinear parabolic equation
	\be
	\left\{
	\begin{array}{l}
	0\=-\Lc u(t,x)-f(x,u(t,x),Du(t,x)\sigma(x)) \;\;,\;(t,x) \in \Oc\x [0,T),  \\
	u|_{\partial_p D}\=g   \;. 
	\end{array}
	\right. \label{PDE}
	\ee 
where we recall that  $\Lc$ denotes the Dynkin operator associated to the diffusion $X$,
	$
	\Lc \psi := \partial_t  \psi +	    \scap{b}{D \psi}  +  \frac12  \Tr{ a D^2\psi   }$ with $a:=\sigma\sigma^*$, 
and 
	$
	\partial_p D :=  ([0,T)\x \partial \Oc) \cup (\{T\}\x \bar \Oc)
	$
is the parabolic boundary of $D$. 

\begin{Proposition}\label{prop u PDE} 
Let \A{HL}, \A{D1},  \A{D2}, \A{C} and \A{Hg} hold. Then the function $u$ has linear growth and is the unique continuous viscosity solution of 
\reff{PDE} in the class of continuous solutions with polynomial growth. 
\end{Proposition} 

A similar result is proved in \cite{DaPa97} but in the elliptic case. For the sake of completeness, we provide a  slightly different complete proof of the viscosity property in the Appendix, where the standard associated comparison result leading to uniqueness is also stated. 
\subsection{Boundary modulus of continuity}\label{sec gradient bound}

 Adapting some barrier techniques for PDEs, we first prove the following bound for the modulus of continuity on the boundary. 
\begin{Proposition}\label{BORNE_GRADIENT_BORD}
Let \A{HL}, \A{D1},  \A{D2}, \A{C} and \A{Hg} hold.  Then, there is  $C_L>0$ such that for all $(t_0,x_0)\in [0,T)\x \partial \Oc $,
\begin{eqnarray}
 \lim_{y\in \cO,\;y\to x_0} \frac{|u(t_0,y)-u(t_0,x_0)|}{\|y-x_0\|} \le C_L.\label{eq borne pour grad}
\end{eqnarray}
 In particular, if the gradient of $u$ exists at $(t_0,x_0)$, it is uniformly bounded.
 \end{Proposition}
 
\textit{Proof.}
Let $(t_0,x_0)\in[0,T)\times \partial\cO $ and  $\cA:=[t_0,T)\times \cN$, where $\cN \subset \cO$ is an open set and $x_0\in \partial \cN $. We only show that, for all $y \in \cN$,
 \begin{eqnarray}
 \frac{u(t_0,y)-u(t_0,x_0)}{\|y-x_0\|}  \le C_L\;.\label{eq borne inf pour grad}
	\end{eqnarray}
The lower bound is obtained similarly. By \A{D2}, there is $\eps>0$ and a family $(e_i)_{i\in \leftB 1,d\rightB}$ such that $x_0+\eps e_i \in \Nc$ for all $i\in \leftB 1,d\rightB$ and span$(e_i,\;i\in \leftB 1,d\rightB)$ $=$ $\R^d$. Thus,  \reff{eq borne pour grad} implies the statement concerning the gradient, whenever it is well defined. We now prove \reff{eq borne inf pour grad}.

\no 1. Assume that there exists a smooth function $\psi: \bar \cA\rightarrow \R$ with first derivative bounded by $C_L$ such that
\begin{enumerate}
\label{bornitude}
\item[(a)] 
 $\psi\ge u$ on $\partial_p \Ac:=([t_0,T)\x \partial \Nc)\cup (\{T\} \x \bar \Nc)$.
\item[(b)] 
\label{neg}
$\Lc \psi(t,x) +f(x,\psi(t,x),D \psi(t,x)\sigma(x)) \le 0 $ for $(t,x)\in \cA$.
\item[(c)] 
\label{eqb}
$\psi(t_0,x_0)=u(t_0,x_0)=g(t_0,x_0)$.
\end{enumerate}
Using Proposition \ref{prop u PDE}  and a standard maximum principle, see Lemma \ref{lem comp edp} in the Appendix, we then derive  that $u\le \psi $ on $\bar \cA$. In view of
(c) this yields
	$$
	\frac{u(t_0,y)-u(t_0,x_0)}{\|y-x_0\|}\le \frac{\psi(t_0,y)-\psi(t_0,x_0)}{\|y-x_0\|}\le C_L\;\;,\;\forall y\in \bar \cN\setminus\{x_0\}\;.
	$$
 
\no 2. It remains to construct a smooth function satisfying (a), (b) and (c). 
Recall that the spatial boundary $\dO $ is compact. Since $u$ is continuous on $\bar D$, see Proposition \ref{prop u PDE}, the compactness assumption \A{D1} ensures the uniform boundedness of $u$ in a neighborhood of $[0,T]\times \dO$.
 
We specify the construction of the barrier function only for $x_0\in \dO\backslash B(\cC,L^{-1})$.  
Indeed, for $x_0\in   B(\cC,L^{-1})$, assumption \A{C} ensures that the diffusion coefficient is uniformly elliptic in a neighborhood of $x_0$. The expression of the barriers below can then be simplified. Namely, we do not need the additional localization with the cone, i.e. we can take $\kappa=0$ in  \reff{eq def psi barriere} below. 

\vs2

Let $y:=y(x_0)$ be the point of $\bar \Oc^c$ associated to $x_0$ by the exterior sphere property, see \A{D2}.  Set $r:=r(x_0)=\|y(x_0)-x_0\|$. Recall that, by assumption,    
	$
	B:=B(y,r) \;\mbox{ satisfies } \; \bar B\cap \bar \Oc=\{x_0\}\;.
	$

It follows from \A{HL} and \A{C}  that 
	\be\label{eq a(x)n(x0)n(x0)}
	\langle a(x) n(x_0), n(x_0)\rangle \ge L^{-1}/2 
	&\mbox{ on the set }& 
	\cD_1:=\{x\in\cO: \|x-x_0\|\le  \eta_L \}\;
	\ee
for some $\eta_L>0$ small enough, but depending only on $L$.

For $x\in\cO$, we now set 
	\b*
	d_B(x):=d(x,\partial B)=\|x-y\|-r
	\e*
 so that $d_B\in C^2(\bar \cO)$ with 
 	\be\label{eq DdB er D2dB}
 	Dd_B(x)=\frac{x-y}{\|x-y\|}\;,\; D^2d_B(x)=\frac{I_d}{\|x-y\|}-\frac{(x-y)^*(x-y)}{\|x-y\|^3} \;
 	\ee
where $I_d$ denotes the identity matrix of $\M^d$.
We now introduce a cone 
	$$
	\cK:=\{x\in\R^d:\ \langle x-y, n(x_0)\rangle \ge \cos(\theta)\|x-y\| \},\ \theta\in[0,\pi/2] 
	$$ 
and 
	$$
	\cD_2:=\{x\in\cO: d_B(x)\le \delta\}\;\;,\;\;\delta >0\;,
	$$
where $\delta\le \delta_L$ small enough to
ensure $\cD_2\subset \cD_1$. 
We finally set
	$
	\cN:=\cO\cap\cK \cap \cD_2\;  
	$ 
and define the barrier function by 
	\be\label{eq def psi barriere}
	\psi(t,x):=g(t,x)+4\alpha (\varphi(x)^{1/2}-\delta^{1/2}) +\kappa \langle x-y, n(x_0) \rangle \left(1-\frac{\langle x-y, n(x_0) \rangle}{\|x-y\|}\right) 
	\ee
for  $(t,x)\in [t_0,T]\times \bar\cN$, where 
	$
	\varphi(x):=\delta+d_B(x)\;.
	$
for some  $(\alpha,\kappa)\in (0,\infty)^2$ to be chosen later on.
\begin{figure}[htb]
\begin{center}
\psfrag{dD}{$\partial O$}
\psfrag{dD1}{$\partial {\cal D}_1$}
\psfrag{dD2}{$\partial {\cal D}_2$}
\psfrag{K}{${\cal K}$}
\psfrag{y}{$y$}
\psfrag{x}{$x$}
\psfrag{n(x)}{$n(x)$}
\psfrag{Theta}{$\theta $}
\includegraphics[
width=0.3\textwidth,angle=0]
{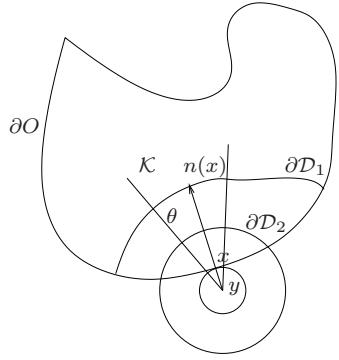}
\caption{Domain for the barrier }
\label{figure:projections}
\end{center}
\end{figure}

2.b.  Since $x_0-y\in$ span$(n(x_0))$, $\psi(t_0,x_0)=u(t_0,x_0)=g(t_0,x_0)$,  so that (c) is
satisfied. 

2.c.  Recall from the beginning of Step 2. that 
	\be\label{eq def M barriere}
	M:=\sup_{(t,x)\in [t_0,T]\times \bar \cD_1}|u(t,x)|\vee \sup_{(t,x) \in [t_0,T]\times \bar \cD_1}|g(t,x)| < \infty\;. 
	\ee
On $\dO\cap \partial\cN$, $\psi(t,x)\ge g(t,x)$. On $\partial \cD_2\cap \partial \cN$,  
	$
	\psi(t,x)\ge -M+4\alpha(2^{1/2}-1)\delta^{1/2}\;. 
	$
Thus, for
\begin{equation}
\label{cond_1}
\alpha\ge \frac{M}{2(2^{1/2}-1)\delta^{1/2}},
\end{equation}
one has $\psi(t,x)\ge u(t,x)$ for $ (t,x)\in [t_0,T]\times  \partial \cD_2 \cap \partial\cN$. 

On $\partial \cK\cap\partial\cN$, we have 	
	$$
	\psi(t,x)\ge -M+\kappa \cos(\theta)\|x-y\|(1-\cos(\theta))\ge -M+\kappa r \cos(\theta)(1-\cos(\theta))\;.
	$$
Hence, for 
\begin{equation}
\label{cond_2}
\kappa\ge \frac{2M}{r\cos(\theta)(1-\cos(\theta))},
\end{equation}
we obtain that $\psi(t,x)\ge u(t,x)$ $\forall (t,x)\in  [t_0,T]\times  \partial\cK\cap
\partial\cN$.  This concludes the proof of (a).

2.d. It remains to show that $\psi$ satisfies (b).  
Set 
	$$
	\Gamma(x):=\langle x-y,n(x_0) \rangle \left(1-\frac{\langle x-y,n(x_0) \rangle}{\|x-y\|} \right) \;,
	$$
and observe that, for some  $C\le C_L$,  
	\be\label{eq borne Gamma barriere}
	\|D \Gamma(x)\|\le C \;,\;\| D^2\Gamma(x)\|\le C /r\;
	\ee
uniformly in $x\in \bar\cN$. 
Define,
\begin{eqnarray*}
 \Theta(t,x)&:=& \Lc\psi(t,x)+f(x,\psi(t,x),D \psi(t,x) \sigma(x)) 
 \\
 &\le& C (1+M+\alpha \varphi(x)^{-1/2}+\kappa(1+r^{-1}))
 -\frac{\alpha}{2}\left\langle a(x)\frac{x-y}{\|x-y\|} ,\frac{x-y}{\|x-y\|}\right\rangle \varphi(x)^{-3/2}
 \\
 &+& \frac{C \;\alpha}{r}\varphi(x)^{-1/2}
 \\
 &\le& C \;(1+M+\kappa(1+r^{-1}))
 -\frac{\alpha}{2}\varphi(x)^{-3/2}\left(\left\langle a(x)\frac{x-y}{\|x-y\|} ,\frac{x-y}{\|x-y\|}\right\rangle-C(1+r^{-1})\varphi(x) \right)\;,
\end{eqnarray*}
recall \A{Hg}, \reff{eq DdB er D2dB}, \reff{eq def M barriere} and \reff{eq borne Gamma barriere}.
For a suitable angle of the cone $\theta$, we shall show below that we can find $\tilde C>0$ such that $\tilde C^{-1}\le C_L$ and 
\begin{equation}
\label{BORNE}
 \left\langle a(x)\frac{x-y}{\|x-y\|} ,\frac{x-y}{\|x-y\|}\right\rangle\ge \tilde C\;\;,\; \forall x\in\bar \cN\;. 
\end{equation}
Recalling that $\varphi(x) \le 2\delta$ for $x\in\bar \cN\subset \Dc_2$,  we get
\begin{eqnarray*}
\Theta(t,x)
\le
C(M+\kappa(1+r^{-1}))-\frac{\alpha}{2}\varphi(x)^{-3/2}\left(\tilde C- 2C(1+r^{-1})\delta \right)\;.
\end{eqnarray*}
For $\delta:=(1/4) \tilde C(C(1+r^{-1}))^{-1}\wedge \delta_L>0$,
we then have 
	$
	\Theta(t,x)\le C(M+\kappa(1+r^{-1}))-\tilde C \alpha 2^{-\frac 72} \delta^{-\frac 32} 
	$.  
It is then clear that $(\alpha,\kappa) $ can be chosen in order to satisfy \reff{cond_1}, \reff{cond_2} and so that $\Theta(t,x)\le 0 $. This shows
(b).  

It remains to prove \reff{BORNE}. This is done by suitably choosing the angle of the cone $\cK$. Let $Z\in \partial B(0,1)$  be such that $Z+y \in \cK$. Introduce
the basis $(n(x_0),(n_i^\bot(x_0))_{i\in\leftB 1,d-1\rightB} )$ where $(n_i^\bot(x_0))_{i\in\leftB 1,d-1\rightB}$ is an orthonormal basis of $\{n(x_0)\}^\bot$ for the euclidean scalar product. Let $(\beta_i)_{i\in \leftB 0,d-1\rightB}$ denote the coefficients of $Z$ in this basis, i.e.  
	$
	Z=\beta_0 n(x_0)+\sum_{i=1}^{d-1}\beta_in_i^\bot(x_0)\;. 
	$
One has, for all $x\in\bar \cN $,
\begin{eqnarray*}
\langle a(x)Z,Z\rangle 
&=&\beta_0^2 \langle a(x)n(x_0),n(x_0)\rangle 
+2\sum_{i=1}^{d-1}\beta_0\beta_i \langle a(x)n(x_0),n_i^\bot(x_0)  \rangle\;
\\
&+&
 \langle a(x)\sum_{i=1}^{d-1}\beta_i n_i^\bot(x_0),\sum_{i=1}^{d-1}\beta_i n_i^\bot(x_0)  \rangle\;
\\
&\ge&\beta_0^2 \langle a(x)n(x_0),n(x_0)\rangle 
+
2\sum_{i=1}^{d-1}\beta_0\beta_i \langle a(x)n(x_0),n_i^\bot(x_0)  \rangle\;.
 \end{eqnarray*}

Since $Z+y\in\cK$ and $\|Z\|=1$, we must have  $\beta_0\ge \cos\theta$, by definition of $\cK$, and therefore $|\beta_i|\le \sin(\theta)$ for all  $i\in\leftB 1,d-1\rightB$.
Hence, \reff{eq a(x)n(x0)n(x0)} and the above equation leads to 
\begin{eqnarray*}
\langle a(x)Z,Z\rangle 
&\ge&
\cos^2(\theta)\frac{L^{-1}}{2}-2(d-1)\sin(\theta)\sup_{x\in \bar \Nc}\|a(x)\|\;,\ \forall x\in \bar \cN.
\end{eqnarray*}
This yields \reff{BORNE} with $\tilde C=\frac{L^{-1}
\cos^2(\theta)}{4} $ for $\theta $ small enough.
\ep


\subsection{Representation and weak regularity of the gradient in the regular uniformly elliptic case}

In the section, we strengthen the initial assumptions and work under: 

\A{D'}:  $\Oc$ is a $C^2$ bounded domain satisfying \A{D1} and \A{D2} for the constant $L$. 

\A{C'}: $a$ is uniformly elliptic with ellipticity constant  $L^{-1}$. 

\A{H'}: the coefficients $b$, $\sigma$, $f$ and $g$ satisfy \A{Hg}-\A{HL} and are uniformly $C^{2}(\bar D)$. 
\\

From now on, given a matrix $M$, we denote by $M^{\cdot j}$ its $j$-th column, viewed as a column vector. 

\begin{Proposition}[Representation of the gradient]\label{REP_GRAD} Let the conditions \A{D'}, \A{C'} and \A{H'} hold. Then, $u \in  C^{0}(\bar D)\cap C^{1,2}(D)$, $Du \in C^0(\bar D)$ and for all $(t,x) \in \bar D$
	\begin{eqnarray}\label{eq_REP} 
	Du(t,x)
	=
	\Esp{
	D u(\tau^{t,x},X^{t,x}_{\tau^{t,x}})\nabla  X^{t,x}_{\tau^{t,x}} V^{t,x}_{\tau^{t,x}}
	+
	\int_{t}^{\tau^{t,x}} \partial_x f(\Theta_s^{t,x}) \nabla X^{t,x}_s V^{t,x}_s ds
	}
	\end{eqnarray}
 where $\nabla  X^{t,x}$ is the first variation process of $X^{t,x}$:
	 $$
	 \nabla  X_s^{t,x}
	 =
	 I_d+\sum_{j=1}^{d}\int_{t}^{s}D \sigma^{\cdot j }(X^{t,x}_v)\nabla X^{t,x}_v dW_v^j+\int_{t}^{s} D b(X^{t,x}_v)\nabla X^{t,x}_v dv\;\;,s\ge t\;,
	 $$ 
and  $V^{t,x}$ is defined by 		
	$$
	V_s^{t,x}
	:=
	\exp\left(
	\int_{t}^{s}\partial_y f(\Theta_v^{t,x}) dv+\int_{t}^{s}\partial_z f(\Theta_v^{t,x}) dW_v
	-\frac12\int_{t}^{s}\|\partial_z f(\Theta_v^{t,x})\|^2dv
	\right)\;  \;\;,\;s\ge t\;,
	$$ 
with $\Theta^{t,x}=(X^{t,x},Y^{t,x},Z^{t,x})$. 
\end{Proposition}

\proof The result is obvious for $(t,x) \in \partial D$. We then assume from now on that $(t,x) \in D$. 
We derive from Theorems 12.16 and 12.10 in \cite{Li05}  and the  definition of Hölder spaces at p. 46 of this reference that $Du \in C^{0}(\bar D)$. 
Let us consider the systems of differential equations obtained by   formally differentiating the PDE \reff{PDE} w.r.t. $(x^i)_{i\in\leftB 1,d\rightB}$. For $i=1,\ldots,d$, we have 
\begin{eqnarray}
0&=&\partial_t v^i+\scap{ b+\sigma^* D_z f(\Theta)+\frac 12 D_{x^i}a^{\cdot i}}{D v^i}  +\frac 12 \Tr{a D^2{v^i}}
\label{PDEvi}\\
&+&\left(D_{x^i}b^i+D_{y} f(\Theta)+\scap{D_z f(\Theta)}{D_{x^i}\sigma^{\cdot i}}\right)v^i+D_{x^i} f(\Theta)+\sum_{k\neq i}h^{i,k}\;,\nonumber\\
&&\mbox{ where } \;\;h^{i,k}=  \left(D_{x^i}b^k+\scap{D_z f(\Theta)}{D_{x^i}\sigma^{\cdot k}}\right) D_{x^k}u+ \sum_{l=1}^d D_{x^i}a^{kl}D_{x^kx^l}u\;\nonumber
\end{eqnarray}
and $\Theta(t,x)=(x,u(t,x),Du\sigma(t,x))$. 

Given $n$ large enough, set $\Oc_n:=\{x+\bar B(0,n^{-1}),\;x\in \Oc^c\}^c$ $\subset \Oc$, $T_n:=T-n^{-1}>0$ and $D_n:=[0,T_n)\x \Oc_n$. 
Note that by construction $\Oc_n$ satisfies a uniform exterior sphere property (with radius $1/2n$).  Then, the PDE \reff{PDEvi} on $D_n$ with the boundary condition $D_{x^i}u$  on $\partial_p D_n=([0,T_n)\x \partial \Oc_n) \cup (\{T_n\}\x \bar \Oc_n)$ admits a unique $C^0(\bar D_n)\cap C^{1,2}(D_n)$ solution $v^i_n$, see Theorem 12.22  in \cite{Li05}. Using the maximum principle, we can then identify $D_{x^i}u$ and $v^i_n$ on $\bar D_n$ by considering the PDE satisfied by   $\eps^{-1}(u(\cdot,x+\eps e_i)- u(\cdot,x))-v^i_n(\cdot,x)$ on $\bar D_n$. Here, $e_i$ is the $i$-th canonical basis vector of $\R^d$, see e.g. Theorem 10 Chapter 3 in \cite{Fr64}. In particular, $Du$ $\in C^0(\bar D_n)\cap C^{1,2}(D_n)$. By a usual localization argument, we then deduce from It\^{o}'s Lemma applied to $Du(\cdot,X^{t,x})\nabla X^{t,x}V^{t,x}$, with $(t,x) \in D_n$, that 
	\b*
	Du(t,x)
	&=&
	\Esp{
	Du(\tau_n, X^{t,x}_{\tau_n})\nabla  X^{t,x}_{\tau_n} V^{t,x}_{\tau_n}
	+
	\int_{t}^{\tau_n} \partial_x f\left(\Theta^{t,x}_s\right) \nabla X^{t,x}_s V^{t,x}_s  ds}
	\;
	\e*
where $\tau_n:=\inf\{s\in [t,T_n]~:~(s,X^{t,x}_s)\notin D_n\}$.  Observe that $\lim_n \tau_n= \tau$ $\Pas$ by continuity  of   $X$.
 We then derive the statement of the Proposition by sending $n\to \infty$, using   the a-priori smoothness of $u$, $Du\in C^{0}(\bar D)$, and the dominated convergence theorem.
\ep

\begin{Remark}{\rm  Note that the various localizations in the previous proof are needed because we do not  assume any compatibility condition on the parabolic boundary, i.e. $\Lc g+f(\cdot,g,\sigma Dg)=0$ on $\partial_p D$. Otherwise, Theorem 12.14 in \cite{Li05} would give  $u\in C^{1,2}(\bar D)$ which would allow to avoid the introduction of the subdomains $\Oc_n$. 

}
\end{Remark}

Observe that, by Proposition \ref{BORNE_GRADIENT_BORD} and the continuity of $Du$ stated in Proposition \ref{REP_GRAD}, we have $\|Du(\tau^{t,x},X^{t,x}_{\tau^{t,x}}) \|\le C_L$. The representation \reff{eq_REP} and standard estimates then give $\|Du\|_{\infty,\bar D}\le C_L$.

\begin{Corollary}\label{coro unif lipschitz} Let \A{D'}, \A{C'} and \A{H'} hold. Then, $\|Du\|_{\infty,\bar D}\le C_L$.
\end{Corollary}

We can now prove Theorem \ref{thm main regu}   under the conditions \A{D'}, \A{C'} and \A{H'}.

\begin{Corollary}\label{coro cas regu} Theorem \ref{thm main regu} holds under the conditions \A{D'}, \A{C'} and \A{H'}.
\end{Corollary}
  
\proof 1. {\bf Proof of \reff{eq main thm regu 3} and \reff{eq main thm regu 2}.} Recalling that $u \in C^{1,2}(D)\cap C^{1}(\bar D)$, see Proposition \ref{REP_GRAD}, we deduce from a standard verification argument that $Z=Du(\cdot,X)\sigma(X)$. Set $(\nabla X, V):=(\nabla X^{0,X_0},V^{0,X_0}) $ and observe that $(\nabla X^{t,X_t}_s, V^{t,X_t}_s)=(\nabla X_s \nabla X_t^{-1} ,V_s V_t^{-1})$ for $s\ge t$, by the flow property. Thus, by Proposition \ref{REP_GRAD},
	\be\label{eq repre Z}
	Z_t
	=
	\Esp{
	Du(\tau, X_{\tau})\nabla  X_{\tau} V_{\tau}
	+
	\int_{t}^{\tau} \partial_x f\left(\Theta_s\right) \nabla X_s V_s  ds~|~\Fc_t}\sigma(X_t)(\nabla X_t V_t)^{-1}\;\;\;,\;t\le \tau\;. 
	\;
	\ee
It then follows from Proposition \ref{BORNE_GRADIENT_BORD} (boundedness of the gradient of $u$), \A{HL} and standard estimates  that $\sup_{t\le \tau} \|Z_t\|$  $\le \xi_L$. This readily implies \reff{eq main thm regu 2}, i.e. $\Esp{\int_\theta^\vartheta \|Z_s\|^p ds~|~\Fc_\theta}$ $\le$ $ \Esp{\xi_L^p|\vartheta-\theta|~|~\Fc_\theta}$, $p=1,2$. By Burkholder-Davis-Gundy's inequality, \A{HL} and Proposition \ref{prop X bar X}, this also yields $\Esp{\sup_{t\in [\theta,\vartheta]}  |Y_t-Y_\theta |^{2p}}$ $\le$  $\Esp{\xi_L^p\;|\vartheta-\theta|^p}$, $p\ge 1$

2. {\bf Proof of \reff{eq regu u}.} By the same arguments as above, we first obtain that $|u(t,x)-u(t,x')|\le C_L |x-x'|$. Moreover, for $t\le t'\le T$, 
	\b*
	u(t,x)-u(t',x)=Y^{t,x}_t-u(t',x)=Y^{t,x}_t-Y^{t,x}_{t'}+ u(t',X^{t,x}_{t'})- u(t',x)\;. 
	\e*
The Lipschitz continuity of $u$ in space (Corollary \ref{coro unif lipschitz}) and standard estimates on SDEs imply that $|\E[u(t',X^{t,x}_{t'})-$ $u(t',x)]|$ $\le$ $C_L|t-t'|^\frac12$. On the other hand,  
	$
	\Esp{|Y^{t,x}_t-Y^{t,x}_{t'}|^2}
	 \le    C_L (t'-t)$,  by the above estimate.   

3. {\bf Proof of \reff{eq main thm regu}.} The bound on $\Rc(Y)_{\Sc^2}^\pi$ follows from  \reff{eq main thm regu 3}. Using \reff{eq repre Z} and  exactly the same arguments as in the proof of Proposition 4.5 in \cite{BoEl05}, see also \cite{MaZh02}, we deduce that 
	$$
 	\sum_{i=0}^{n-1} \Esp{\int_\ti^\tip \|Z_{t}-Z_{\ti}\|^2dt}
 	\le
	C_L\;h\;\;,
	$$
which implies 
	$
 	\sum_{i=0}^{n-1} \Esp{\int_\ti^\tip \|Z_{t}-\hat Z_{\ti}\|^2dt}
 \le
	C_L\;h\;\;
	$
since $\hat Z$ is the best approximation of $Z$ in $L^2(\Omega\x [0,T])$  by an element of $\Hc^2$ which is constant on each time interval $[\ti,\tip)$.  \ep

\subsection{Regularization procedure: proof of Theorem \ref{thm main regu}  in the general case}\label{subsec regularization procedure}

\textbf{Step 1. Truncation of the domain:} We first prove that Theorem \ref{thm main regu} holds under the conditions \A{D1}, \A{D2}, \A{C'} and \A{H'}. 

Let $\phi $ be a $C^{\infty}$   density function with compact support on $\R^d$. Given $\eps>0$, we define $\Delta_{\eps}:=\eps^{-d}\phi(\eps^{-1}\cdot)\star(d\wedge d_{\eps^{-1}})^+$ where $d_{\eps^{-1}}$ denotes the algebraic distance to $\partial B(X_0,\eps^{-1})$ and $\star$ denotes the convolution.  Set $\Oc_{\eps}:=\{x\in \R^d~:~\Delta_{\eps}(x)>0\}$ and $D_{\eps}:=[0,T)\x \Oc_{\eps}$. It follows from the compact boundary assumption that  $\partial \Oc\subset \bar \Oc_\eps$, for $\eps$ small enough. Note that $\Oc_\eps$ is bounded, even if $\Oc$ is not. Let $(Y^\eps,Z^\eps)$ be defined as in 
\reff{eq def BSDE intro} with $\Oc_\eps$ in place of $\Oc$ and $\tau^\eps$ be the first exit time of $(\cdot,X)$ from $D_\eps$. 
Observe  that, by continuity of $X$, $\tau^\eps\to \tau$ $\Pas$ Since,  by \A{Hg}, \A{HL} and Theorem 1.5 in \cite{Pa98},  
	\b*
	\|Y-Y^\eps\|^2_{\Sc^2}+\|Z-Z^\eps\|^2_{\Hc^2} &\le&  
	C_L\Esp{|g(\tau,X_\tau)-g({\tau^\eps},X_{\tau^\eps})|^2+\int_{\tau\wedge {\tau^\eps}}^{\tau\vee {\tau^\eps}} f(X_s,Y_s,Z_s)^2 ds}
	\\
	&\le& C_L \Esp{ \int_{\tau\wedge {\tau^\eps}}^{\tau\vee {\tau^\eps}} (1+\|X_s\|^2+  |Y_s|^2+ \|Z_s\|^2 )ds }\;,
	\e*
we deduce from Proposition \ref{prop X bar X} and a dominated convergence argument  that $\|Y-Y^\eps\|^2_{\Sc^2}+\|Z-Z^\eps\|^2_{\Hc^2}$ $\to$ $0$. Since the domain $\Oc_\eps$ satisfies \A{D'}, we can apply Corollary \ref{coro cas regu} to $(Y^\eps,Z^\eps)$. Recalling that the associated constants depend only on $L$ and are uniform  in $\eps$, we thus obtain the required controls on $(Y,Z)$.  
Let $u^\eps$ be the solution of \reff{PDE} associated to $D_\eps$. The above stability result, applied to general initial conditions, implies that $u^\eps \to u$ pointwise on $\bar D$.  Corollary \ref{coro cas regu} thus implies that $u$ satisfies \reff{eq regu u}. 

\vs2

\textbf{Step 2. Regularization of the coefficients:} We now prove that Theorem \ref{thm main regu} holds under the conditions \A{D1}, \A{D2}, \A{C}, \A{HL} and \A{Hg}. 

For  $\varepsilon>0$, define $b_\eps$, $\sigma_\eps$ and $f_\eps$  by
	\b*
	(b_\eps, \sigma_\eps,f_\eps)(x,y,z):=(b,\sigma,f)\star \varepsilon^{-2d+1}\phi(\varepsilon^{-1}(x,y,z)) 
	\e*
where  $\phi $ is a $C^{\infty}$  density function with compact support on $\R^d\x\R\x\R^d$.
Let us consider the FBSDE
	\be
	\label{BSDE_epsi}
	\begin{cases}
	X_t^\varepsilon=x+\int_{0}^{t}b_\eps(X_s^\varepsilon)ds+\int_{0}^{t}\sigma_\eps(\Xe_s)dW_s+\sqrt{\varepsilon}\tilde W_t\;,
	\\
	\Ye_t=g(\te,\Xe_{\te})+\int_{t\wedge \te}^{\te}f_\eps(\Xe_s,\Ye_s,\Ze_s)ds-\int_{t\wedge\te}^{\te}\Ze_sdW_s-\int_{t\wedge \te}^{\te}\tilde Z_s^\eps d\tilde W_s\;,
	\end{cases}
	\end{eqnarray}
where $(\tilde W_t)_{t\ge 0}$ is an additional $d$-dimensional Brownian motion independent of $W$ and 
	$$
	\te:=\inf \{s\ge 0: (s,\Xe_s)\not\in D   \} \;.
	$$ 
This system satisfies the conditions of Step 1. Therefore, the estimates of Theorem \ref{thm main regu} can be applied to $(Y^\eps,Z^\eps)$. Note that the associated constant depends only on $L$ and are uniform  in $\eps$. Moreover, it follows  from \A{HL} and  Theorem 1.5 in \cite{Pa98} that 
	\b*
	\|Y-Y^\eps\|^2_{\Sc^2}+\|Z-Z^\eps\|^2_{\Hc^2} &\le&  C_L\Esp{ |g(\tau,X_\tau)-g(\tau^\eps,X^\eps_{\tau^\eps}) |^2 + \int_0^T \| X_s-X^\eps_s\|^2ds }\\
	&+&  \Esp{\int_{\tau \wedge \tau^\eps}^{\tau \vee \tau^\eps}  (|f(X_s,Y_s,Z_s)|+|f_\eps(X^\eps_s,Y_s,Z_s)|)^2ds }+L\;\eps\;. 
	\e*
Clearly, $X^\eps\to X$ in $\Sc^2$. Since $f$ and $g$ are Lipschitz continuous, $f$ and $f_\eps$ have linear growth and $(X,X^\eps,Y,Z)$ is bounded in  $\Sc^2\x\Sc^2\x\Sc^2\x\Hc^2$, it suffices to check that $\tau^\eps\to \tau$ in probability to obtain the required controls on $(Y,Z)$. This is implied by the non-characteristic boundary condition of \A{C}, see e.g. the proof of Proposition 3 in \cite{GoMe04}. The control  \reff{eq regu u} is obtained by arguing as above. 
\ep
 
  
\appendix\section*{Appendix: Proof of Proposition \ref{prop u PDE}}

In the following, we use the notations
	\b*
	u^*(t,x)=\limsup_{(s,y)\in D \to (t,x)} u(s,y) \;,\;u_*(t,x)=\liminf_{(s,y)\in D \to (t,x)} u(s,y)\;\;,\;(t,x)\in \bar D.
	\e*

The statement of Proposition \ref{prop u PDE} is a direct consequence of Lemmas \ref{lem A1} 	and \ref{lem comp edp} below.

\begin{Lemma}\label{lem A1}Let the conditions of Proposition \ref{prop u PDE} hold.  Then, the function $u$ has linear growth and  $u^*$ (resp. $u_*$) is a viscosity subsolution (resp. supersolution) of \reff{PDE} with the terminal conditions $u^*\le g$ (resp. $u_*\ge g$) on $\partial_p D$. 
\end{Lemma}

\proof 1. The linear growth property property is an immediate consequence of Proposition \ref{prop X bar X}.
 
2. It remains to prove that $u^*$ and $u_*$ are respectively sub- and supersolution of \reff{PDE} with the boundary conditions $u^*\le g$ and $u_*\ge g$ on $\partial_p D$.  We concentrate on the supersolution property, the subsolution property would be derived similarly. The proof is standard, as usual  we argue by contradiction. 
 Let $(t_0,x_0) \in [0,T]\x \bar \Oc$ and $\vp \in C^2_b$ be such that $0=\min_{(t,x)\in \bar D} (u_*-\vp)(t,x)=(u_*-\vp)(t_0,x_0)$ where the minimum is assumed, w.l.o.g., to be strict on $\bar D$. 
Assume that  
 \b*
 \left(-\Lc \vp(t_0,x_0)-f(x_0,\vp(t_0,x_0),D\vp\sigma(t_0,x_0))\right)\1_{(t_0,x_0)\in D} + (\vp-g)(t_0,x_0)\1_{(t_0,x_0)\in \partial_p D}=:-2\zeta<0 \;.~~ 
 \e*
Recall from \A{D2} that  if $x_0\in \partial \Oc$ then we can find an open ball $B_0\subset \Oc^c$ such that $\bar B_0\cap \bar \Oc=\{x_0\}$. If $x_0 \in \partial \Oc$, we denote by $d_{B_0}$  the algebraic distance to $B_0$.
On $\bar D$, we set 
	\b*
	\tilde \vp(t,x) &=& \vp(t,x)-(\sqrt{T-t})\1_{t_0=T} -d(x)\left(1-\frac{d(x)}{\eta}\right) \1_{x_0\in \partial \Oc \setminus B(\Cc,L^{-1}) }
										 \\
&&-d_{B_0}(x)\left(1-\frac{d_{B_0}(x)}{\eta}\right) \1_{x_0\in \partial \Oc\cap B(\cC,L^{-1}) }  \;,
	\e*
for some $\eta>0$. Observe that $(t_0,x_0)$ is still a strict minimum of $(u_*-\tilde \vp)$ on $V_\eta\cap  \bar D$ for some open neighborhood $V_\eta$ of $(t_0,x_0)$ on which $(d_{B_0}\vee d)\le \eta/2$ if $x_0 \in\partial \Oc $. Without loss of generality, we can then assume that 
	\be\label{eq sur sol contra u - vp}
	u\ge u_*\ge \tilde \vp +\zeta \;\mbox{ on } \; \partial V_\eta\setminus \bar D^c\;,
	\ee 
while 
	\be\label{eq sur sol contra loc bis}
	\tilde \vp \le \vp\le g -\zeta \;\mbox{ on } \; \bar V_\eta\cap  \partial_p  D \;,\;\mbox{ if } (t_0,x_0)\in \partial_p D\;. 
	\ee 
Moreover, observe that for $F$ equal to $d$ or $d_{B_0}$, $D(F(1-F/\eta))=DF(1-2\eta^{-1} F)$ and $D^2(F(1-F/\eta))=(1-2\eta^{-1}F)D^2F-2\eta^{-1}DF^* DF$ where $\|DF\|=1$. Thus, \A{C} implies that, for $\eta$ and  $V_\eta$ small enough,  
 \be\label{eq sur sol contra loc}
  -\Lc \tilde \vp-f(\cdot,\tilde \vp,D\tilde \vp\sigma)  \le - \zeta<0 \; \mbox{ on } \; V_\eta\cap \bar D \;.
 \ee
 Let $(t_n,x_n)_n$ be a sequence in $D\cap V_\eta$ such that $(t_n,x_n,u(t_n,x_n))\to (t_0,x_0,u_*(t_0,x_0))$. Let $(X^n,Y^n,Z^n)$ be the solution of \reff{eq def SDE intro}-\reff{eq def BSDE intro} associated to the initial conditions $(t_n,x_n)$ and define $\theta_n$ as the first exit time of $D\cap V_\eta$  by $(\cdot,X^n)$. By applying Itô's Lemma on $\tilde \vp$ and using  \reff{eq sur sol contra loc bis}, \reff{eq sur sol contra loc},  \reff{eq sur sol contra u - vp} and   the identity $u=g$ on $\partial_p D$, we get
	\b*
	\tilde \vp(t_n,x_n)
	&=&
	-\chi+ 
	u(\theta_n,X^n_{\theta_n})  +  \int_{t_n}^{\theta_n } (f(X^n_s,\tilde \vp(s,X^n_s),D\tilde \vp\sigma(s,X^n_s))-\eta_s) ds 
	\\
	&&-\; \int_{t_n}^{\theta_n }  D\tilde \vp\sigma(s,X^n_s)dW_s\;,
	\e*
where $\chi$ is a   bounded random variable satisfying $\chi\ge \zeta$ $\Pas$ and $\eta$ is an adapted process in $L^2$ such that $\eta\ge \zeta$ $dt\x d\P$-a.e.
Following the standard argument of the proof of Theorem 1.6 in  \cite{Pa98}, we deduce that $\tilde \vp(t_n,x_n)\le Y^{t_n,x_n}_{t_n}-\zeta e^{-LT}=u(t_n,x_n)-\zeta e^{-LT}$.  Since $\tilde \vp(t_n,x_n)-u(t_n,x_n)\to 0$, this leads to a contradiction. \ep\\

We now state a comparison theorem for the PDE \reff{PDE}. The proof is quite standard, see e.g. \cite{CrIsLi92}, but we give it for the sake of completeness. 

\begin{Lemma}\label{lem comp edp} Let the conditions of Proposition \ref{prop u PDE} hold.   Fix $t_0 \in [0,T)$ and $\Nc\subset \Oc$ an open set.  Let $U$ (resp. $V$) be an upper-semicontinuous subsolution (resp. lower-semicontinuous supersolution) with polynomial growth of \reff{PDE} on $\Ac:=[t_0,T)\x \Nc$  such that $V\ge U$ on $\partial_p \Ac:=([t_0,T)\x \partial \Nc)\cup (\{T\}\x \bar \Nc)$. Then,   $V\ge U$ on $\bar \Ac$.
\end{Lemma} 

\proof Fix $\rho>0$ and observe $\tilde U$ and $\tilde V$ defined by $\tilde U(t,x)=U(t,x)e^{\rho t}$ and $\tilde V(t,x)=V(t,x)e^{\rho t}$ are sub- and supersolution of 
	\be
	0\=\rho\psi(t,x) -\Lc \psi(t,x)-e^{\rho t} f(x,e^{-\rho t}\psi(t,x),e^{-\rho t}D\psi(t,x)\sigma(x)) \;\;,\;(t,x) \in [t_0,T)\x\Nc\;.  \label{PDE tilde}
	\ee 
As usual we argue by contradiction and assume that $\sup_{(t,x)\in A} (\tilde U(t,x)-\tilde V(t,x))>0$. 
Define 
	\b*
	\beta(t,x):= e^{- \kappa t} (1+\|x\|^{2p}) \;\;,\;(t,x)\in \bar \Ac
	\e*
for $p\in \N^*$ such that $(|U(t,x)|+|V(t,x)|)/(1+\|x\|^{p})$ is bounded on $\bar \Ac$, and $\kappa>0$ to be chosen later on.  
For all $\eps>0$ small enough, we can then find   $(t_\eps,x_\eps) \in \bar \Ac$ such that 
	 \be\label{eq def txeps}
	 \sup_{(t,x)\in \Ac} (\tilde U(t,x)-\tilde V(t,x)-2\eps \beta(t,x))=:(\tilde U(t_\eps,x_\eps)-\tilde V(t_\eps,x_\eps)-2\eps \beta(t_\eps,x_\eps))>0\;.
	 \ee
Clearly, $(t_\eps,x_\eps)\notin \partial_p \Ac$ since $\tilde U\le \tilde V$ on $\partial_p \Ac$. For $n\in \N^*$, let $(t_n,x_n,y_n) \in [t_0,T]\x \bar \Nc^2$ be a maximum point of 
	\b*
	(\tilde U(t,x)-\tilde V(t,y)-\eps (\beta(t,x)+\beta(t,y)) - \left(|t-t_\eps|^2+\|x-x_\eps\|^4+n\|x-y\|^2\right)\;.  
	\e*
It is easy to check, see e.g. Proposition 3.7 in \cite{CrIsLi92}, that 
	\begin{equation}
	\tilde U(t_n,x_n)-\tilde V(t_n,y_n)\to (\tilde U-\tilde V)(t_\eps,x_\eps)
	 \And 
	|t_n-t_\eps|^2+\|x_n-x_\eps\|^4+n\|x_n-y_n\|^2 \to  0\;.
	\label{eq visco lim n}
	\end{equation}
Since $(t_\eps,x_\eps)\in \Ac$, we can assume that  $(t_n,x_n)\in \Ac$ for all $n\in \N^*$, after possibly passing to a subsequence. It then follows from Ishii's Lemma, Theorem 8.3 in \cite{CrIsLi92}, that  we can find 	 real coefficients $a_{n}$,
$b_{n}$ and symmetric matrices $\Xc_n$ and
$\Yc_n$ such that
    \b*
    \left(a_{n},p_n, \Xc_n\right)
    \;\in\; \bar {\Pc}^+_{\bar \Nc}\tilde U(t_n,x_n) &\And&
    \left(b_{n},q_n, \Yc_n\right)
    \;\in\; \bar {\Pc}^-_{\bar \Nc}\tilde V(t_n,y_n)\;,
    \e*
see \cite{CrIsLi92} for the standard notations $\bar {\Pc}^+_{\bar \Nc}$ and $\bar {\Pc}^-_{\bar \Nc}$,
where
    \b*
    p_n := 2n (x_n-y_n)+4(x_n-x_\eps)\|x_n-x_\eps\|^2+ \eps D\beta(t_n,x_n)
&,&
    q_n :=  2n (x_n-y_n)- \eps D\beta(t_n,y_n)
    \e*
and 
    \begin{equation}\label{ineq b and A}
    a_{n} - b_{n} = 2 (t_n - t_\eps) +\eps  \left(\partial_t\beta(t_n,x_n)+\partial_t\beta(t_n,y_n)\right)
   \;,\;
    \left(
    \begin{array}{c c}
    \Xc_n&0\\ 0&-\Yc_n
    \end{array}
    \right)
    \le   A_n  + n^{-3} (A_n)^2~~~~
    \end{equation}
with
    \b*
    A_n  &:=& 2n\left(\begin{array}{cc} I_d &-I_d\\
-I_d &I_d \end{array}\right)+\eps	\left(
    \begin{array}{c c}
     D^2\beta(t_n,x_n)  &  0\\
     0  &        D^2\beta(t_n,y_n) 
    \end{array}
    \right)\\
&&
     +\left(\begin{array}{cc}4 I_d\|x_n-x_\eps\|^2+8(x_n-x_\eps)^*(x_n-x_\eps)  &0 \\0 & 0 \end{array} \right),
    \e*
where $I_d$ is the identity matrix of $\M^d$. 
Since $\tilde U$ and $\tilde V$ are sub- and supersolution of \reff{PDE tilde}, it follows that 
	\b*
	\rho\left(\tilde U(t_n,x_n)-\tilde V(t_n,y_n)\right)
	&\le&
	a_n-b_n + \scap{b(x_n)}{p_n}-\scap{b(y_n)}{q_n}  +  \frac12  \Tr{ a(x_n) \Xc_n  - a(y_n) \Yc_n   }
	\\
	&+&
	\left(f(x_n,U(t_n,x_n),e^{-\rho t_n}p_n\sigma(x_n))
	-
	f(y_n,V(t_n,y_n),e^{-\rho t_n}q_n\sigma(y_n))
	\right)e^{\rho t}\;. 
	\e* 
We then deduce from \A{HL}, \reff{ineq b and A}, \reff{eq visco lim n}, and standard computations that
	\b*
	\rho\left(\tilde U(t_n,x_n)-\tilde V(t_n,y_n)\right)
	&\le &
	L \left(\tilde U(t_n,x_n)-\tilde V(t_n,y_n)\right)
	+
	2\eps  \left(\Lc \beta(t_\eps,x_\eps) + L \| \sigma(x_\eps)D\beta \|\right)
	+
	o_n(1)
	\;.
	\e*
Taking $\rho>2L$ and $\kappa$ large enough so that $ \Lc \beta+ L \| \sigma D\beta \| \le -\frac \kappa 2\exp(-\kappa T)$ on $\bar \Ac$, which is possible thanks to \A{HL}, we finally obtain 
	\b*
	\frac12 \rho\left(\tilde U(t_n,x_n)-\tilde V(t_n,y_n)\right)
	&\le & -	\kappa \exp(-\kappa T)\eps 
	+
	o_n(1)\;, 
	\e* 
which contradicts \reff{eq def txeps} for $n$ large enough, recall \reff{eq visco lim n}.

\ep



{\small
\begin{thebibliography}{aa12}
 
 
\bibitem{BaPa02} Bally  V., and G. Pages (2002). A quantization algorithm
for solving discrete time multidimensional optimal stopping
problems. {\sl Bernoulli}, 9 (6), 1003-1049.
 
\bibitem{BoCh06} Bouchard B. and J.-F. Chassagneux (2006). Discrete time approximation for continuously and discretely reflected BSDE's. Preprint LPMA, University Paris 6.

\bibitem{BoEl05} Bouchard  B. and R. Elie  (2005). Discrete time approximation of decoupled Forward-Backward SDE with jumps. To appear in {\sl Stochastic Processes and their Applications}.

\bibitem{BoTo04} Bouchard  B. and N. Touzi  (2004).  Discrete-Time Approximation and Monte-Carlo Simulation of Backward Stochastic Differential Equations.
{\sl Stochastic Processes and their Applications}, 111 (2), 175-206.
  
 
 
\bibitem{CrIsLi92} { Crandall M. G., H. Ishii and P.-L. Lions} (1992). {
User's guide to viscosity solutions of second order Partial
Differential Equations}. {\sl Amer. Math. Soc.}, 27, 1-67.


\bibitem{DaPa97} Darling R. W. R.  and E. Pardoux (1997). {BSDE} with random terminal time. {\sl Annals of Probability}, 25(3), 1135-1159.

\bibitem{DeMe06}   Delarue F. and S. Menozzi (2006).  A forward backward algorithm for quasi-linear PDEs. {\sl Annals of Applied Probability},  16 (1), 140-184.

\bibitem{DeMe07}   Delarue F. and S. Menozzi (2007).   An interpolated Stochastic Algorithm for Quasi-Linear PDEs. {\sl Mathematics of Computation}, 261-77 (2008), 125-158. .

 
\bibitem{Fr85} Freidlin M. (1985). {\sl Functional integration and partial differential equations}. Annals of Mathematics Studies, Princeton University Press. 


\bibitem{Fr64}  Friedman A.  (1964). Partial Differential Equations
of Parabolic Type. {Prentice Hall}.

\bibitem{GiTr98}   Gilbarg D.  and N. S. Trudinger (1998). 
{\sl Elliptic Partial Differential Equations of Second Order}. 
Springer-Verlag, Berlin.
 
\bibitem{Go98} Gobet E. (1998). {\sl Schéma d'Euler pour diffusions tuées. Application  aux options barrière}. Phd Thesis, University Paris VII.

\bibitem{Go00} Gobet E. (2000). Weak approximation of killed diffusion using Euler schemes. {\sl Stochastic Processes and their Applications}, 87, 167-197.

\bibitem{GoLa07} Gobet E. and C. Labart (2007).  Error expansion for the discretization of Backward Stochastic Differential Equations. {\sl Stochastic Processes and Applications}, 117 (7), 803-829.

\bibitem{GoLeWa05} Gobet  E., J.P. Lemor   and  X. Warin  (2006).
Rate of convergence of empirical regression method for solving
generalized BSDE. {\sl Bernoulli},  12 (5), 889-916. 

\bibitem{GoMe07} Gobet E. and S. Menozzi (2007). Stopped diffusion processes: overshoots and boundary correction. Preprint PMA, University Paris 7. 

\bibitem{GoMe04} Gobet E. and S. Menozzi (2004). Exact approximation rate of killed hypoelliptic diffusions using the discrete Euler scheme. {\sl Stochastic Processes and their Applications}, 114 (2), 201-223.

\bibitem{KaSh90} Karatzas I. et S.E. Shreve (1990). {\sl Brownian motion and stochastic calculus}. Springer Verlag.

\bibitem{Ko00}   Kobylanski M. (2000). {B}ackward {S}tochastic {D}ifferential {E}quations and {P}artial {D}ifferential {E}quations with quadratic growth. {\sl 
Annals of Probability}, 28(2), 558-602.

\bibitem{LaSoUr68} Ladyzenskaja O.A.,  V.A. Solonnikov and N.N. Ural'ceva (1968). {\sl Linear and quasi-linear equations of parabolic type}. Trans. Math. Monog., vol. 23, AMS, Providence. 


\bibitem{Li05}  Lieberman G. M. (2005). Second Order Parabolic Differential Equations. {World Scientific}.


\bibitem{LoSc01}
 Longstaff F. A. and  R. S. Schwartz  (2001). Valuing American
Options By Simulation~: A simple Least-Square Approach. {\sl
Review of Financial Studies},  {\  14}, 113-147.
 





\bibitem{MaZh02}  Ma J. and Zhang  J. (2002). Path Regularity of Solutions to Backward
Stochastic Differential Equations. {\sl Probability Theory and
Related Fields}, 122, 163-190.
 
 
\bibitem{MaZh05}  Ma J. and Zhang  J. (2005). Representations and regularities for solutions
to BSDEs with reflections. {\sl Stochastic Processes and their Applications}, 115, 539-569.

\bibitem{MiTr01} Milstein  G.N. and M.V. Tretyakov   (2001). Numerical solution of Dirichlet problems for nonlinear parabolic equations by a probabilistic approach. {\sl IMA J. Num. Anal.},  21 (4), 887-917.


\bibitem{Pa98} Pardoux  E. (1998). Backward stochastic differential equations and  viscosity solutions of semilinear parabolic and elliptic PDE's of second order. In {\sl Stochastic Analysis and Related Topics: The Geilo Workshop 1996}. L. Decreusefond, J. Gjerd, B. Oksendal, and A.S. Ustünel (eds.), Birkhäuser,   79-127. 
 
\bibitem{PaPe92} Pardoux  E. and S. Peng  (1992). Backward stochastic differential
equations and quasilinear parabolic partial differential equations. {\sl Lecture Notes in Control and Inform. Sci},
176, 200-217.
 
\bibitem{Pe91} Peng S. (1991).  Probabilistic interpretation for systems of quasilinear parabolic partial differential equations. {\sl Stochastics and Stochastics reports}, 37, 61-74.
 
  
 
\bibitem{Zh04} Zhang J. (2004).   A numerical scheme for BSDEs. {\sl Annals of
Applied Probability},  14 (1), 459-488.
\end{thebibliography}
}
\end{document}